\documentclass[11pt,a4paper,openbib]{article}
\usepackage[margin=26mm,marginparwidth=21mm,marginparsep=2.5mm]{geometry}
\usepackage[T1]{fontenc}
\usepackage[utf8]{inputenc}
\usepackage{lmodern,microtype,amsmath,amssymb,amsthm,mathtools,mathrsfs,bm}
\usepackage{enumitem,xcolor,longtable,hyperref}
\hypersetup{colorlinks=true,linkcolor=blue,citecolor=blue,urlcolor=blue,pdftitle={Local scaling limit of the four-dimensional uniform spanning tree},pdfauthor={Runsheng Liu},pdfkeywords={uniform spanning tree, loop-erased random walk, Brownian Kesten tree, scaling limit}}
\setlist{leftmargin=1.8em,itemsep=0.2em,topsep=0.3em}
\numberwithin{equation}{section}
\newtheorem{theorem}{Theorem}[section]
\newtheorem{proposition}[theorem]{Proposition}
\newtheorem{lemma}[theorem]{Lemma}
\newtheorem{corollary}[theorem]{Corollary}
\theoremstyle{definition}
\theoremstyle{remark}\newtheorem{remark}[theorem]{Remark}
\newcommand{\Pp}{\mathbb P}\newcommand{\Ee}{\mathbb E}

\newcommand{\R}{\mathbb R}\newcommand{\Z}{\mathbb Z}

\newcommand{\one}{\mathbf 1}\newcommand{\dd}{\,\mathrm d}
\newcommand{\Law}{\operatorname{Law}}

\title{Local scaling limit of the four-dimensional uniform spanning tree}
\author{Runsheng Liu\thanks{School of Mathematical Sciences, Peking University. \href{mailto:liurunsheng@pku.edu.cn}{liurunsheng@pku.edu.cn}}}
\begin{document}\maketitle
\begin{abstract}
In this paper, we prove that the uniform spanning tree of $\mathbb Z^4$, seen from the origin, converges locally to an infinite spatial Brownian Kesten tree: a Brownian spine onto which independent spatial Brownian trees are grafted at a constant rate. The convergence holds jointly for the intrinsic metric, the volume measure and the Euclidean positions of the vertices, after rescaling these quantities by suitable constant multiples of $n^{-1}$, $n^{-2}(\log n)^{1/3}$ and $n^{-1/2}(\log n)^{-1/6}$, respectively, under the local rooted spatial Gromov-weak topology. The rescaled space-time occupation measure converges to that of infinite canonical super-Brownian motion. We also identify the limiting law of the rescaled intrinsic-ball volume and obtain sharp first-moment asymptotics.
\end{abstract}

{\small
\noindent\textit{2020 Mathematics Subject Classification.} Primary 60K35; Secondary 60G50, 60J80.\par
\noindent\textit{Keywords.} Uniform spanning tree, loop-erased random walk, Brownian Kesten tree, scaling limit.\par
}

\section{Introduction}\label{min:setup}

\subsection[Motivation and background]{Motivation and background}\label{intro:background}

This paper studies the local scaling limit of the uniform spanning tree (UST) of $\mathbb Z^4$, recording its intrinsic metric, volume measure and spatial embedding jointly. Wilson's algorithm \cite{Wilson1996} constructs UST branches from loop-erased random walks (LERWs), so their geometry is governed by random-walk intersections. To add a vertex, the walk is stopped on hitting the existing tree and its loop-erasure is attached there. Four is the critical dimension for these intersections \cite{Lawler1991}, and logarithmic corrections enter the relevant scales. The infinite-volume spanning forest is connected in dimensions $d\le4$ and has infinitely many components in dimensions $d\ge5$ \cite{Pemantle1991}. In four dimensions, it can be constructed by Wilson's algorithm rooted at infinity \cite{BenjaminiLyonsPeresSchramm2001}.

Previous work has identified logarithmic corrections to intrinsic distance and volume estimates for four-dimensional USTs. On the torus $(\mathbb Z/m\mathbb Z)^4$, Schweinsberg \cite{Schweinsberg2009} proved convergence of distance arrays for finitely many independent uniform vertices to those of the Brownian continuum random tree (CRT) \cite{Aldous1991}, with distance scale $m^2(\log m)^{1/6}$ up to deterministic factors bounded above and away from zero. On the infinite lattice, Hutchcroft and Sousi \cite{HutchcroftSousi2023} determined logarithmic corrections to the intrinsic-radius and volume tails of the past of a vertex. Halberstam and Hutchcroft \cite{HalberstamHutchcroft2024} established intrinsic-ball volume growth of order $n^2(\log n)^{-1/3+o(1)}$ and associated geometric and spectral corrections. These results leave open the local scaling limit of the infinite tree around a vertex. We identify this limit as the spatial Brownian Kesten tree in the local rooted spatial Gromov-weak topology, with distances, volume and positions rescaled together. Our proof also gives sharp first-moment asymptotics for intrinsic-ball volume.

The Brownian Kesten tree is suggested by the mean-field picture: weak interactions between separated branches should produce a critical branching structure, with logarithmic corrections in dimension four. A one-ended critical tree viewed from a fixed vertex should resemble Kesten's tree \cite{Kesten1986}, a critical Galton-Watson tree conditioned to survive, with an infinite spine and independent critical trees attached along it. Nachmias and Peres \cite{NachmiasPeres2022} established the Poisson(1) Kesten tree as the local limit of USTs on finite simple regular graphs of diverging degree, rooted at a uniform vertex. The Brownian Kesten tree arises by rescaling finite-variance critical Galton-Watson trees conditioned to survive \cite{Duquesne2009}. In the continuum picture, the spine and ancestral paths carry Brownian spatial labels, while volume is carried by the trees grafted along the spine. Van der Hofstad \cite[Sections~6.2--6.3]{VanDerHofstad2006} conjectured convergence of the rescaled USF $r$-point functions to the moment measures of infinite canonical super-Brownian motion (ICSBM) above four dimensions and raised the corresponding question in the critical dimension. Corollary~\ref{min:volume-occupation-law} below gives an ICSBM limit for the space-time occupation measure of the four-dimensional UST.

In two dimensions, the conformally invariant UST and LERW limits of \cite{Schramm2000,LawlerSchrammWerner2004} were complemented by measured spatial-tree tightness \cite{BarlowCroydonKumagai2017} and convergence through the Peano curve and its intrinsic contour functions \cite{HoldenSun2018}. In three dimensions, convergence in a local Gromov-Hausdorff-Prokhorov topology recording the embedding was established along dyadic meshes in \cite{AngelCroydonHernandezTorresShiraishi2021}. The subsequent LERW results in \cite{HernandezTorresLiShiraishi2024} extend this convergence to arbitrary vanishing mesh sequences.

For a class of rapidly mixing finite graphs, including discrete tori in dimensions $d\ge5$, finite sampled-distance convergence to the CRT was proved in \cite{PeresRevelle2005} and strengthened to Gromov-Hausdorff-Prokhorov convergence in \cite{ArcherNachmiasShalev2024}. For uniform spanning forests on $\mathbb Z^d$, $d\ge5$, upper and lower tail bounds for intrinsic-ball volumes on the scale $n^2$ were obtained in \cite{BarlowJarai2019}, while \cite{Hutchcroft2020} established mean-field geometric and spectral exponents for a broad class of graphs.

Our argument uses four-dimensional LERW estimates to identify both the Brownian branches and their attachment rate. The relation between walk steps and LERW length and the invariance principle in \cite[Chapter~7]{Lawler1991}, together with the logarithmic correction in \cite{Lawler1995}, describe the branch scale. The renormalized escape laws and stationary two-sided LERW of \cite{LawlerSunWu2019}, and the capacity law and finite-segment escape limits of \cite{Markering2026}, identify capacity on the scale $n/(\log n)^{2/3}$. We show that local escape probabilities yield the same limiting attachment rate. The capacity coefficient thus becomes the grafting rate per unit intrinsic length, linking the potential theory of LERW to the branching structure of the limit.

\subsection{Model and main result}\label{intro:main-result}

Let $\mathcal U$ be the wired uniform spanning forest of $\mathbb Z^4$. It is almost surely a connected, one-ended tree \cite{BenjaminiLyonsPeresSchramm2001}. Root it at $0$, write $d_{\mathcal U}$ for graph distance, and denote the unique ray from $x$ towards infinity by $\gamma^x=(\gamma^x_j)_{j\ge0}$, with $\gamma^x_0=x$. By Wilson's algorithm rooted at infinity \cite{BenjaminiLyonsPeresSchramm2001}, the ray $\gamma^0$ is an infinite LERW from $0$: the chronological loop-erasure of a simple random walk run forever, which is well defined because the walk is transient. We call $\gamma^0$ the \emph{spine}. Unless otherwise stated, all random walks below are nearest-neighbour simple random walks on $\mathbb Z^4$, so their one-step covariance is $\Sigma=I_4/4$. Let $\mathsf P$ be the law of the infinite loop-erasure of an unconditioned simple random walk, and let $\alpha=\mathsf E_{\mathsf P}X_\infty\in(0,\infty)$ be the mean limiting normalized escape probability defined in \eqref{prelim:raw-escape-normalization}. It determines the relation between the number of random-walk steps and the resulting LERW length.

We use the scales
\begin{equation}\label{min:scales}
L_n=(\log(n+2))^{1/3},\qquad N_n=\frac{nL_n}{\alpha},\qquad R_n=\sqrt{N_n},\qquad b_n=\frac{\alpha L_n}{n^2},
\end{equation}
and study the rooted spatial measured tree
\begin{equation}\label{min:rescaled-tree}
\mathscr U_n= \left(\mathcal U,n^{-1}d_{\mathcal U}, b_n\sum_{x\in\mathbb Z^4}\delta_x, \phi_n,0\right),\qquad \phi_n(x)=R_n^{-1}x.
\end{equation}
Edges may be replaced by intervals of length $1/n$, with linearly interpolated labels.

\emph{The capacity constant.} For a finite set $A\subset\mathbb Z^4$, our capacity convention is
\begin{equation}\label{min:capacity-convention}
\operatorname{cap}(A)=\sum_{x\in A}\Pp_x(\tau_A^+=\infty), \qquad \tau_A^+=\inf\{j\ge1:W_j\in A\}.
\end{equation}
There is no factor of the vertex degree in this definition. The positive constant $\kappa$ is defined by the stationary escape-law expression \eqref{min:kappa}. Corollary~\ref{min:capacity-prop} proves that, for ordinary infinite LERW $\eta$ from $0$,
\begin{equation}\label{min:kappa-normalization}
\frac{L_n^2}{n}\operatorname{cap} (\eta[0,\lfloor nt\rfloor])\longrightarrow\kappa t \quad\text{in probability},\qquad t>0.
\end{equation}
 In words, a LERW of $n$ edges has capacity about $\kappa n/(\log n)^{2/3}$, so each of its vertices carries an escape probability of order $(\log n)^{-2/3}$ on average. Neither $\alpha$ nor $\kappa$ is known explicitly, but their product is: $\alpha\kappa=\pi^2/12$ by \eqref{min:escape-constant-check}.

\emph{The topology in words.} We want to compare $\mathscr U_n$ with a continuum tree while keeping track of tree distances, volume and spatial positions at the same time. The simplest way to do this is to sample. Pick $k$ points according to the volume measure, which is not normalized, and record all their tree distances, to each other and to the root, together with their spatial positions. Two measured spatial trees are close if, for every $k$, the resulting measures on these finite arrays are close. This is a rooted, spatially marked version of the Gromov-weak topology. Since $\mathscr U_n$ is infinite, we apply it to intrinsic balls around the root.

The limiting object is the spatial Brownian Kesten tree $\mathcal K_{\kappa,\Sigma}$: an infinite ray, called its spine, labelled by Brownian motion, with independent spatial Brownian-tree excursions attached along it. The precise definitions of the topology and the limit tree are given below, after the two statements.

\begin{theorem}[Local spatial limit]\label{min:main}
With the scales in \eqref{min:scales},
\[
\mathscr U_n\Longrightarrow\mathcal K_{\kappa,I_4/4}
\]
in local rooted spatial Gromov-weak topology.
\end{theorem}

The theorem also identifies the space-time occupation measure. Let $\mathcal X=(\mathcal X_t)_{t\ge0}$ be ICSBM issued from the origin, with $\mathcal X_0=0$, branching mechanism $\psi(q)=\kappa q^2/2$ and spatial covariance $I_4/4$. Its construction is recalled in Section~\ref{prelim:continuum}.

\begin{corollary}[ICSBM occupation measure]\label{min:volume-occupation-law}
For every $T>0$,
\begin{equation}\label{min:occupation-limit}
b_n\sum_{x\in B_{\mathcal U}(0,\lfloor nT\rfloor)}\delta_{(d_{\mathcal U}(0,x)/n,\,x/R_n)}\Longrightarrow\mathbf1_{[0,T]}(t)\,dt\,\mathcal X_t(dx)
\end{equation}
in distribution for the weak topology on finite measures on $[0,T]\times\mathbb R^4$. In particular, $b_n|B_{\mathcal U}(0,\lfloor nT\rfloor)|\Longrightarrow V_T:=\int_0^T\mathcal X_t(\mathbb R^4)\,dt$, whose Laplace transform is given in Proposition~\ref{prelim:icsbm-volume}.
\end{corollary}

The precise definitions are as follows. We use the rooted, finite-measure spatially marked version of the marked Gromov-weak topology \cite[Definitions~2.3--2.4 and Theorems~1--2]{DepperschmidtGrevenPfaffelhuber2011}, with localization by intrinsic balls as in \cite[Lemma~2.6 and Definition~2.7]{AthreyaLohrWinter2016}. A rooted finite spatial metric-measure space is a tuple $\mathsf X=(X,d,m,\phi,o)$, where $(X,d)$ is complete and separable, $m$ is a finite Borel measure, $o\in X$, and $\phi:X\to\mathbb R^4$ is a Borel label map. All spaces considered here have $\phi(o)=0$. We regard these labels as the finite marked measure $(\mathrm{id},\phi)_*m$ on $X\times\mathbb R^4$, in the general marked-measure state space of \cite[Section~2.1]{FoutelRodier2024}.  We identify spaces by rooted, measure-preserving isometries of $\operatorname{supp}(m)\cup\{o\}$ that preserve labels $m$-almost everywhere and at the root. For $k\ge1$, put $x_0=o$ and define its rooted marked distance-matrix measure by
\begin{equation}\label{min:distance-matrix}
\nu_k^{\mathsf X}= \left[ (x_1,\ldots,x_k)\longmapsto \big((d(x_i,x_j))_{0\le i<j\le k}, (\phi(x_i))_{0\le i\le k}\big) \right]_*m^{\otimes k}.
\end{equation}
Rooted marked Gromov-weak convergence means weak convergence of \eqref{min:distance-matrix} for every $k$.  For spaces whose measure is finite on every bounded intrinsic ball, we use the corresponding local topology, obtained by restricting $m$ to intrinsic balls about $o$. Here and below a restricted ball retains the root, its inherited metric and labels.  For finite spaces, let $d_{\mathrm{r,mGP}}$ be the rooted marked Gromov-Prokhorov metric: take the infimum, over common isometric embeddings, of the maximum of the distance between the roots and the Prokhorov distance between the embedded marked measures. It metrizes rooted marked Gromov-weak convergence. Writing $B_R(\mathsf X)=\mathsf X\!\restriction_R$ for restriction to the closed root ball, the local topology is metrized by
\begin{equation}\label{min:local-metric}
d_{\mathrm{loc}}(\mathsf X,\mathsf Y) =\int_0^\infty e^{-R} \bigl(1\wedge d_{\mathrm{r,mGP}} (\mathsf X\!\restriction_R,\mathsf Y\!\restriction_R)\bigr)\,dR.
\end{equation}
For random spaces, the measures $\Ee[\nu_k^{\mathsf X}]$ are their rooted marked moment measures.

To construct $\mathcal K_{\kappa,\Sigma}$, let $\xi$ be a Brownian motion with covariance $\Sigma$, and assign the label $\xi_s$ to the point of the spine at distance $s$ from the root. Spatial Brownian trees are grafted onto the spine according to a Poisson point process with intensity $\kappa\,ds\,\mathbb N_{\xi_s}(d\mathcal E)$, where $\mathbb N_x$ is the excursion measure of a spatial Brownian tree started at $x$. All the volume lives on the grafted trees. The same constant $\kappa$ appears as the grafting rate here and as the capacity per unit length in \eqref{min:kappa-normalization}. Section~\ref{prelim:continuum} gives its construction and normalization.

\begin{remark}\label{intro:remarks}
(a) In concrete terms: fix $T>0$ and $k\ge1$, and sample $k$ vertices from the ball $B_{\mathcal U}(0,nT)$ according to the counting measure weighted by $b_n$. The resulting random measure on their tree distances, divided by $n$, and their positions, divided by $R_n$ (that is, the measure $\nu_k$ of \eqref{min:distance-matrix} for the rescaled ball), converges in law to the corresponding random measure for the ball of radius $T$ in $\mathcal K_{\kappa,I_4/4}$.

(b) The proof also provides a sharp asymptotic estimate for the first moment of the ball volume. Proposition~\ref{bridge:initial-kernel} implies
\[
\mathbb E\,|B_{\mathcal U}(0,\lfloor nT\rfloor)|\sim\frac{\pi^2T^2}{24\alpha^2}\,n^2(\log n)^{-1/3}\qquad(n\to\infty),
\]
which matches the expected volume $\kappa T^2/2$ of the root ball of radius $T$ in the limit tree; see \eqref{min:kesten-first-moment}. Corollary~\ref{min:volume-occupation-law} gives the volume convergence, and Proposition~\ref{prelim:icsbm-volume} identifies the limiting law.

(c) The Gromov-weak topology records only the measured geometry. It does not impose Hausdorff convergence of branches that carry little or no mass. Upgrading Theorem~\ref{min:main} to a local Gromov-Hausdorff-Prokhorov statement, as obtained in three dimensions in \cite{AngelCroydonHernandezTorresShiraishi2021}, would require additional tightness estimates that we do not address here.

(d) Convergence of the untruncated higher moment measures would additionally require uniform integrability of the corresponding powers of the rescaled ball volume. The deletion estimate controls only the expected removed volume.
\end{remark}

\emph{Choice of scales.} In four dimensions, a random walk of $N$ steps leaves a loop-erasure of only about $\alpha N/(\log N)^{1/3}$ edges. Indeed, step $k$ of the walk is the last visit to a vertex of the loop-erasure exactly when the rest of the walk never returns to the loop-erased path built up to time $k$. By Lawler, Sun and Wu \cite{LawlerSunWu2019}, this has probability about $\alpha(\log k)^{-1/3}$; see \eqref{prelim:raw-escape-normalization} and Lemma~\ref{prelim:calibrated-clock}. Hence about $N_n=nL_n/\alpha$ walk steps are needed to produce $n$ LERW edges, and the resulting branch has the diffusive spatial size $R_n=\sqrt{N_n}$ of the walk that produced it. The volume factor $b_n$ is the inverse of the typical volume of an intrinsic ball of radius $n$. In particular, $b_nN_n=L_n^2/n$. A random walk run for about $N_n$ steps produces a loop-erasure of about $n$ edges. Dividing the walk step count by $N_n$ and the LERW edge count by $n$ therefore gives the same limiting parameter, with spatial covariance $I_4/4$. The four rescalings have different roles:
\begin{center}
\small
\begin{tabular}{@{}llll@{}}
Quantity & Original unit & Rescaled value & Order of the unit at scale $n$\\\hline
Tree distance & number of tree edges & $n^{-1}d_{\mathcal U}$ & $n$\\
Walk duration & number of walk steps & $k/N_n$ & $N_n\asymp n(\log n)^{1/3}$\\
Spatial position & lattice vertex & $R_n^{-1}x$ & $R_n\asymp n^{1/2}(\log n)^{1/6}$\\
Volume of a vertex set & number of vertices & $b_n|A|$ & $b_n^{-1}\asymp n^2(\log n)^{-1/3}$
\end{tabular}
\end{center}
\emph{Why $b_n$.} By Wilson's algorithm, the ball $B_{\mathcal U}(0,n)$ consists of pieces of branches grafted onto the first $n$ vertices of the spine. A random walk started at a spine vertex escapes the nearby part of the spine with probability of order $(\log n)^{-2/3}=L_n^{-2}$, roughly one factor $L_n^{-1}$ for each of the two halves of the spine around its starting point. The occupation identity of Corollary~\ref{min:wilson-occupation} sums contributions over attachment points and surviving walk times, with a time range of order $N_n$ corresponding to LERW length of order $n$. The expected volume of the ball is therefore of order $n\cdot L_n^{-2}\cdot N_n=n^2/(\alpha L_n)=b_n^{-1}$. This agrees with the growth $n^2(\log n)^{-1/3+o(1)}$ found by Halberstam and Hutchcroft, and Remark~\ref{intro:remarks}(b) provides the corresponding sharp first-moment asymptotic.

\emph{Capacity and the grafting rate.} The constant governing the limit is a capacity per unit length, as suggested by the following potential-theoretic picture. For a finite set, capacity controls the probability that a distant random walk hits it, and the conditional hitting distribution is approximately the equilibrium measure divided by the capacity. This suggests that branches coming from far away should be grafted onto a long segment of the spine at a rate governed by its capacity.

\begin{remark}
For $d\ge5$, the same strategy should yield a local scaling limit for the component of the origin in the uniform spanning forest on $\Z^d$. More precisely, letting $\mathcal T_d$ be the tree containing $0$, and choosing $b_n^{(d)}=\alpha_dn^{-2}$, $N_n^{(d)}=\alpha_d^{-1}n$ and $R_n^{(d)}=\alpha_d^{-1/2}n^{1/2}$, we similarly define
\[\mathscr T_n^{(d)}= \left(\mathcal T_d,n^{-1}d_{\mathcal T_d}, b_n^{(d)}\sum_{x\in V(\mathcal T_d)}\delta_x, \phi_n^{(d)},0\right),\qquad \phi_n^{(d)}(x)=(R_n^{(d)})^{-1}x.\]
The expected convergence is
\[\mathscr T_n^{(d)}\Longrightarrow\mathcal K_{\kappa_d,I_d/d}\]
under the same topology as in Theorem~\ref{min:main}, with spatial labels in $\R^d$. The constants $\alpha_d,\kappa_d$ above are determined by
\[\alpha_d:=\mathbb P(W[1,\infty)\cap\eta[0,\infty)=\varnothing)\]
and
\[\lim_{m\to\infty}\frac{\operatorname{cap}(\eta[0,m])}{m}=\kappa_d,\quad\mathrm{a.s.},\]
where $W$ is a SRW and $\eta$ is an independent infinite LERW on $\Z^d$, both started from $0$. Thanks to \cite[Theorem 7.7.2]{Lawler1991} and \cite[Theorem 1.4]{Markering2026}, both constants are positive.

As in the four-dimensional argument, one would first prove Campbell-moment convergence for retained vertex measures and then control the removed mass to obtain convergence of the full measured tree. However, the convergence of fixed-level $r$-point functions would require additional estimates.
\end{remark}

\subsection{Proof strategy}\label{intro:strategy}

The proof uses the method of moments for metric measure spaces (Lemma~\ref{prelim:marked-moments}), applied to a restriction of the vertex measure whose removed mass tends to zero in mean. For each $k$, we study the \emph{$k$-point Campbell measure} of the ball of radius $T$. This is defined as follows: integrate over ordered $k$-tuples of vertices in the ball against the $k$-fold product of the volume measure, record the tree distances, spatial positions and ancestral paths of the sampled vertices, and take expectation. The volume measure is not normalized, so the total mass of the $k$-point Campbell measure is the $k$th moment of the normalized ball volume. The Brownian Kesten tree satisfies an exact recursion in $k$: sampling one more vertex amounts to attaching one more Brownian branch to the exposed spine and branches, at rate $\kappa$ per unit length (Corollary~\ref{min:continuum-campbell}). We prove that the retained discrete Campbell measures obey the same recursion in the limit. The steps are as follows.

\emph{Step 1: Wilson insertion for a volume sample.} Fix an intrinsic radius $T$.  Wilson's algorithm and loop-erasure reversal express the contribution of one new vertex as a sum over attachment points and random-walk durations. The walk starts at the attachment and is killed when it returns to the exposed tree. We record both the number of edges in its loop-erasure and the entire loop-erased path. This permits us to add samples successively without losing their genealogical information.

\emph{Step 2: identification of the attachment rate.} The normalization of one insertion can already be seen without its path labels. Let $p_j$ be the probability that a random walk started at the spine vertex $\eta_j$ avoids the nearby piece $\eta[(j-c_n)_+,j+c_n]$ of the spine during steps $1,\ldots,r_n$. The cutoffs $c_n\ll r_n\ll n$, defined in \eqref{not:cutoffs}, are smaller than $n$ by powers of $L_n$. For $0<\varepsilon<M$, the measure with these weights on attachment positions and rescaled walk durations factors exactly as
\[
b_n\sum_{j=0}^{\lfloor nT\rfloor}\sum_{\varepsilon\le k/N_n\le M}p_j\,\delta_{j/n}\otimes\delta_{k/N_n} =\left(\frac{L_n^2}{n}\sum_{j=0}^{\lfloor nT\rfloor}p_j\delta_{j/n}\right)\otimes\left(\frac1{N_n}\sum_{\varepsilon\le k/N_n\le M}\delta_{k/N_n}\right).
\]
The first factor converges to $\kappa\,ds$ by Proposition~\ref{bridge:row-convergence}, and the second is a Riemann sum for $dt$. This explains both the coefficient $\kappa$ and the integration over the new branch length. The estimates below justify using $p_j$ in place of the probability of avoiding the whole path and retain the entire attached path. The restriction to a ball becomes $s+t\le T$.

We first determine the weights of the attachment points. Capacity sums escape probabilities from a whole finite LERW segment, whereas the insertion comparison uses a short segment around each attachment. We show that the two sums have negligible difference after normalization. The capacity law and stationary two-sided escape laws then give the limiting attachment measure $\kappa\,\mathrm ds$.

\emph{Step 3: the Brownian limit of the new branch.} The early local avoidance event gives the leading survival weight, while cut-time comparisons asymptotically separate it from the later increments that produce the Brownian branch. To identify the attached path as well, we replace avoidance of the whole path by an event involving only a short initial segment of the testing walk. Cut times separate this initial segment from the later loop-erasure. The invariance principle, with the loop-erasure parametrized by its edge count, supplies a Brownian branch. Occupation estimates control very short and very long testing-walk durations. This proves Proposition~\ref{bridge:initial-kernel}, including its conditional expected ball-volume limit.

\emph{Step 4: moment bounds for retained samples.} Proposition~\ref{bridge:initial-kernel} controls conditional expected ball volume in $L^p$ for $p<3/2$, which does not supply all the moments needed for iteration. Fix $K>2\kappa T^2$ and retain a vertex when the conditional expected normalized volume of its radius-$2T$ ball, given its full ray towards infinity, is at most $K$; impose the same condition at the root. The next insertion measure then has mass at most $(k+1)K$, and the retained $k$-point Campbell measure has mass at most $k!K^k$. The conditional volume limit and mass transport show that the removed normalized volume tends to zero in mean (Section~\ref{a-truncation-with-vanishing-deleted-mass}).

\emph{Step 5: insertion along a finite history.} After finitely many samples, record their exposed paths and genealogy as a history. Mass transport and the preceding insertion bounds show that the centered full future ray of every retained sample has marginal measure bounded by $k!K^k$ times ordinary LERW law (Lemma~\ref{closure:bounded-rows}). Each added branch, read backwards, is a prefix of its endpoint's future ray. Its next-insertion contribution, with avoidance required only of that ray, is therefore a restriction of the one-ray kernel from Proposition~\ref{bridge:initial-kernel}; reflecting the attachment coordinate and translating by the attachment position give the required Brownian limit. The rest of the exposed tree outside this full ray is a finite set of $O(n)$ vertices. Near recorded endpoints its effect is controlled by the small attachment mass; away from them, spatial separation and a weighted hitting estimate make its effect vanish. The separation follows because four-dimensional Brownian motion has no double points \cite{DvoretzkyErdosKakutani1950}, so the finitely many limiting Brownian arcs are spatially separated away from recorded endpoints. This proves the one-step limit using only convergence of the existing history.

\emph{Step 6: Campbell convergence and the local limit.} Proposition~\ref{closure:campbell-convergence} proves convergence of these restricted Campbell measures by induction on $k$. Lemma~\ref{prelim:marked-moments}, together with the factorial volume moments of the limiting tree, identifies the law of the restricted measured ball. Restoring the removed measure and then letting the radius increase proves Theorem~\ref{min:main}.

\subsection{Organization of the paper}\label{intro:organization}

In Section~\ref{sec:preliminaries}, we introduce the notation and collect preliminary results on LERW and the Brownian Kesten tree. Section~\ref{min:capacity-section} proves convergence of the local escape measure in Proposition~\ref{bridge:row-convergence}. Section~\ref{sec:bridge} proves the initial-spine insertion limit and the conditional ball-volume asymptotic in Proposition~\ref{bridge:initial-kernel}. Section~\ref{sec:closure} introduces the restricted Campbell measures, transfers the one-ray estimates to every sampled future ray, controls hits on the other exposed branches, and completes the induction in Proposition~\ref{closure:campbell-convergence}. Its final argument deduces Theorem~\ref{min:main} from that proposition and the marked metric-measure moment criterion.

\medskip

\noindent {\bf Acknowledgements:} We thank Yifan Gao, Xinyi Li and Satomi Watanabe for useful discussions and helpful comments on this work. This work is supported by the National Key R\&D Program of China (No.~2021YFA1002700) and Beijing Natural Science Foundation (JQ26001).

\medskip

\noindent {\bf AI Statement:} Large language models were used in proving technical lemmas and improving writing quality. The author takes full responsibility for the correctness of this paper.

\section{Preliminaries}\label{sec:preliminaries}\label{min:identities}

Section~\ref{prelim:notation} introduces our notation and provides an index of frequently used symbols. Section~\ref{prelim:lerw} gives the Wilson occupation identity, relates walk steps to LERW length, and records the Brownian path limit and bounds on local escape-probability sums and occupation. Section~\ref{prelim:continuum} constructs the spatial Brownian Kesten tree, identifies its ICSBM occupation measure and ball-volume law, and presents its volume-sampling recursion and the rooted marked moment criterion.

\subsection{Notation}\label{prelim:notation}

We collect general conventions here. The LERW parametrization, escape probabilities and sampling measures are defined where they are first used.

\paragraph{Constants and comparisons.}
The letter $C\ge1$ denotes a finite deterministic constant whose value may change between occurrences. Subscripts, as in $C_{T,p}$, indicate dependence on the listed fixed parameters. Such constants are independent of the scaling index $n$, spatial locations and realizations of paths or exposed trees unless explicitly stated otherwise.  Named constants, such as $C_{\mathrm{time}}$, keep their values once chosen. We write $f_n=O_{\theta}(g_n)$ when $|f_n|\le C_{\theta}g_n$, $f_n=o(g_n)$ when $f_n/g_n\to0$, and $f_n\asymp g_n$ when their ratio is bounded above and away from zero. Parameters are fixed in limits as $n\to\infty$ unless a different order of limits is specified.

\paragraph{Sets and paths.}
We use $a\vee b=\max\{a,b\}$, $a\wedge b=\min\{a,b\}$ and $a_+=a\vee0$. The notation $|x|$ is the Euclidean norm, while $|A|$ is the cardinality of a finite set. The closed ball in a metric space $(X,d)$ is $B_X(x,R)=\{y:d(x,y)\le R\}$. We omit the subscript when the metric is specified by the context. For an indexed path $\omega$, $\omega[a,b]$ denotes its segment and also its range in intersection events. Its reversal is $\overleftarrow\omega$, its number of edges is $|\omega|$, and $\oplus$ denotes concatenation with the common endpoint counted once.

\paragraph{Probability and finite measures.}
We write $\mathbb P,\mathbb E$ for probability and expectation, $\operatorname{Law}(X)$ for a law, and $\mathbf1_A$ for an indicator. For a random walk $X=(X_j)_{j\ge0}$, its natural filtration is $\mathcal F_m^X:=\sigma(X_0,\ldots,X_m)$, $m\ge0$. For a finite nonnegative measure $\mu$ on $E$, write $\mu(F)=\int F\,d\mu$ and denote its total mass by
\[
\|\mu\|:=\mu(E)=\int_E1\,\mu(dx).
\]
For a kernel $Q$, $\|Q(x,\cdot)\|$ denotes the mass of its measure at $x$. Write $f_*\mu$ for pushforward by $f$ and $\Longrightarrow$ for weak convergence of finite measures. For finite measures on a metric space, set
\[
\|\mu-\nu\|_{\rm TV} =\sup_{\|F\|_\infty\le1}|\mu(F)-\nu(F)|,\qquad d_{\rm BL}(\mu,\nu) =\sup_{\substack{\|F\|_\infty\le1\\\operatorname{Lip}(F)\le1}} |\mu(F)-\nu(F)|.
\]
The first supremum is over measurable tests and the second over Lipschitz tests. For a nonnegative measure, $\|\mu\|=\|\mu\|_{\rm TV}$.  The underlying metric is specified with each path or history space below.

\paragraph{Index of frequently used notation.}
The following symbols are used across several sections. The last column says where each is defined.
{\small
\setlength{\LTpre}{0.4em}\setlength{\LTpost}{0.4em}
\begin{longtable}{@{}p{0.24\textwidth}p{0.54\textwidth}p{0.15\textwidth}@{}}
\textbf{Symbol} & \textbf{Meaning} & \textbf{Defined in}\\\hline\endhead
$L_n,\ N_n,\ R_n,\ b_n$ & $(\log(n+2))^{1/3}$; walk steps for $n$ LERW edges; spatial scale $\sqrt{N_n}$; mass per vertex. Note $b_nN_n=L_n^2/n$ & \eqref{min:scales}\\
$\Sigma=I_4/4$ & covariance of one step of simple random walk & \S\ref{intro:main-result}\\
$q_n$ & $\lfloor nT\rfloor$, the intrinsic radius in lattice units & \S\ref{prelim:lerw}\\
$e_A(x)$, $\operatorname{cap}(A)$ & escape probability $\Pp_x(\tau_A^+=\infty)$; capacity $\sum_{x\in A}e_A(x)$ & \eqref{min:capacity-convention}\\
$\mathsf P$, $\eta$ & law of the infinite LERW from $0$; a path with this law & \S\ref{intro:main-result}\\
$X_m$, $X_\infty$, $\alpha$ & normalized escape probability of $\eta$, its limit, and $\alpha=\mathsf E_{\mathsf P}X_\infty$ & \eqref{prelim:raw-escape-normalization}\\
$\mathsf Q$, $\widehat\eta$ & escape-biased law $d\mathsf Q=(X_\infty/\alpha)\,d\mathsf P$; stationary two-sided LERW & \eqref{min:density}\\
$X_\infty^+$, $\kappa$ & escape limit for $\eta[1,\infty)$; capacity per unit length $\kappa=\mathsf E_{\mathsf Q}[X_\infty X_\infty^+]$, with $\alpha\kappa=\pi^2/12$ & \eqref{min:kappa}\\
$\sigma_j$, $\rho_m$ & walk-step index of the $j$th LERW vertex; inverse index $\rho_m=\max\{j:\sigma_j\le m\}$ & \S\ref{prelim:lerw}\\
$c_n\ll r_n\ll u_n$ & window radius (LERW edges); early-avoidance horizon and decoupling time (walk steps) & \eqref{not:cutoffs}\\
$\mathcal C_j^\eta$, $e_j$, $p_j$ & window $\eta[(j-c_n)_+,j+c_n]$; probability that a walk from $\eta_j$ avoids it for all positive times, resp.\ during steps $1,\ldots,r_n$ & \S\ref{prelim:lerw}\\
$Z_{n,T}$, $E_{n,T}$ & $(L_n^2/n)\sum_{j\le q_n}e_j$ and $(L_n^2/n)\sum_{j\le q_n}p_j$ & \S\ref{prelim:lerw}\\
$Q_{n,\mathcal S,T}$, $\overline Q_{n,\mathcal S,T}$ & discrete insertion measures for an exposed tree $\mathcal S$ (the second with atoms at exposed vertices) & \eqref{min:wilson-row}\\
$\mathbb N_x$, $\mathcal K_{\kappa,\Sigma}$ & spatial Brownian-tree excursion measure; Brownian Kesten tree & \S\ref{prelim:continuum}\\
$Q_{\mathcal T,T}$, $\mathfrak C_{k,T}$ & continuum insertion measure; continuum $k$-point Campbell measure & \eqref{min:continuum-row}, Cor.~\ref{min:continuum-campbell}\\
$I^\eta_{n,T}$ & insertion measure on the initial spine & Prop.~\ref{bridge:initial-kernel}\\
$Y_{n,t}(\gamma^x)$ & conditional expected normalized volume of $B_{\mathcal U}(x,\lfloor nt\rfloor)$ given $\gamma^x$ & Prop.~\ref{bridge:initial-kernel}, \S\ref{sec:closure}\\
$g_n$, $G_n$ & retention indicator $\mathbf1\{Y_{n,2T}\le K\}$ and its value at $0$ & \S\ref{sec:closure}\\
$V_n$, $A_n$, $D_n$ & normalized volume of the root ball; its retained and removed parts & \S\ref{sec:closure}\\
$\mathfrak C^{\,g}_{n,k,T}$ & retained discrete $k$-point Campbell measure & \eqref{closure:eq10_1}\\
$H$, $\mathcal S(H)$, $\widehat H_n$ & lattice history, its exposed tree, its rescaling $\operatorname{Scale}_n(H)$ & \S\ref{sec:closure}\\
$R^i$, $\Gamma^i$ & full future ray from sample $x_i$; that ray translated to start at $0$ & \S\ref{closure:ray-reference}\\
$Q_n^{\mathrm{ray}}$ & insertion measure requiring avoidance only of the assigned sample's full ray & \eqref{closure:reference-kernel}
\end{longtable}}

\subsection{Loop-erased random walk and Wilson occupation}
\label{prelim:lerw}\label{sec:clocks}

Loop-erasure reversibility turns vertex sums into occupation by a walk from the attachment point. We establish the Brownian limit in the LERW length parameter, moment bounds for local escape-probability sums, and a bound on occupation at large walk-step indices.

Throughout this subsection, $S$ and $W$ are nearest-neighbour simple random walks on $\mathbb Z^4$, and $B$ is Brownian motion with covariance $\Sigma=I_4/4$. The subscripts $x$ on probability and expectation specify the starting point; $\mathbb E_x^\Sigma$ denotes Brownian expectation. Write $\operatorname{LE}(\omega)$ for chronological loop-erasure. For a vertex set $A$ put
\[
\tau_A=\inf\{k\ge0:W_k\in A\},\qquad \tau_A^+=\inf\{k\ge1:W_k\in A\},\qquad e_A(x)=\Pp_x(\tau_A^+=\infty).
\]
Thus $e_A(x)$ is the probability that the walk avoids $A$ at every positive time. Capacity sums these escape probabilities as in \eqref{min:capacity-convention}. The heat kernel and Green functions are
\[
p_k(x,y)=\Pp_x(W_k=y),\qquad G_D(x,y)=\Ee_x\sum_{k<\tau_{D^{\mathrm c}}}\one_{\{W_k=y\}}, \qquad G=G_{\mathbb Z^4}.
\]

\paragraph{Stopped spatial paths.}\label{not:stopped-paths}
We record a branch by its rescaled number of edges and its spatial trajectory. A stopped path is a pair $(t,\beta)$ with lifetime $t$, where $\beta:[0,\infty)\to\mathbb R^4$ is continuous and constant after $t$. Convergence means convergence of lifetimes and locally uniform convergence of trajectories. The bounded metrics we use are
\[
\begin{aligned}
d_{\rm path}(\beta,\gamma) &=\sum_{a=1}^\infty2^{-a} \left(1\wedge\sup_{0\le v\le a}|\beta(v)-\gamma(v)|\right),\\
d_{\rm stop}((t,\beta),(s,\gamma)) &=\tfrac12\left[(1\wedge|t-s|)+d_{\rm path}(\beta,\gamma)\right].
\end{aligned}
\]
Spatial translations are isometries. If we also record attachment coordinates or random-walk times, we use the maximum of the bounded coordinate metrics. The corresponding $d_{\rm BL}$ is the metric defined in Section~\ref{prelim:notation}.

\paragraph{Wilson occupation with the branch path recorded.}

The following occupation identity is a corollary of classical LERW reversibility \cite[Theorem~4.6]{HutchcroftSousi2023}; see also \cite[Section~3, Propositions~3.1--3.5]{Lawler2017}.

\begin{corollary}[Occupation identity with the loop-erased path recorded]\label{min:wilson-occupation}
Let $\mathcal S\subset\mathbb Z^4$ be a deterministic set and $a\in \mathcal S$. Let $X$ be simple random walk from $z\notin \mathcal S$, with $\tau_{\mathcal S}=\inf\{j\ge0:X_j\in \mathcal S\}$, and let $W$ start at $a$, with $\tau_{\mathcal S}^+=\inf\{j\ge1:W_j\in \mathcal S\}$. For every nonnegative function $F$ on finite oriented simple paths,
\begin{equation}\label{min:occupation}
\sum_{z\notin \mathcal S}\Ee_z\!\left[ \mathbf1_{\{\tau_{\mathcal S}<\infty,\,X_{\tau_{\mathcal S}}=a\}} F\big(\overleftarrow{\operatorname{LE}(X[0,\tau_{\mathcal S}])}\big) \right]= \Ee_a\!\left[\sum_{j=1}^{\tau_{\mathcal S}^+-1} F\big(\operatorname{LE}(W[0,j])\big)\right].
\end{equation}
Both sides may equal $+\infty$. The same identity holds conditionally on a random $\mathcal S$, with fresh test walks independent of the information used to expose $\mathcal S$.
\end{corollary}

In words: the left-hand side runs Wilson's algorithm from every vertex $z\notin\mathcal S$ and collects the branches that attach to $\mathcal S$ at $a$, while the right-hand side runs a single walk from the attachment point $a$ and records its loop-erasure at every time before it returns to $\mathcal S$. The right-hand side is easier to analyse because all walks start from one known point. Its time parameter $j$ plays the role of the duration of the branch-generating walk.

We apply this identity to a connected tree $\mathcal S$ already exposed by Wilson's algorithm, containing $0$. Such a tree may consist of the initial infinite ray and finitely many added branches. Write $h_{\mathcal S}(w)=d_{\mathcal S}(0,w)$ for the number of tree edges from $0$ to $w$. For a test walk started at $w\in\mathcal S$, let
\[
\lambda_{w,j}=\operatorname{LE}(W[0,j]),\qquad J_j=|\lambda_{w,j}|,
\]
where $\lambda_{w,j}$ is recorded as an ordered sequence of lattice vertices and $J_j$ counts its edges. For $T>0$, put $q_n=\lfloor nT\rfloor$ and define
\begin{align}\label{min:wilson-row}
Q_{n,\mathcal S,T}(F) &=b_n\sum_{\substack{w\in \mathcal S\\h_{\mathcal S}(w)\le q_n}} \sum_{j\ge1}\Ee_w\!\left[ \mathbf1_{\{j<\tau_{\mathcal S}^+,\,h_{\mathcal S}(w)+J_j\le q_n\}} F(w,\lambda_{w,j})\right],\notag\\
\overline Q_{n,\mathcal S,T}(F) &=Q_{n,\mathcal S,T}(F) +b_n\sum_{\substack{w\in \mathcal S\\h_{\mathcal S}(w)\le q_n}}F(w,\mathbf0_w).
\end{align}
Here $\mathbf0_w=(w)$ is the singleton lattice path, with zero edges. Both kernels record the attachment vertex and finite path in lattice units. Their rescaled images divide intrinsic coordinates by $n$ and spatial labels by $R_n$. We call these insertion measures. They count vertices with weight $b_n$. For every nonnegative $F$, Wilson's algorithm and \eqref{min:occupation} identify $\overline Q_{n,\mathcal S,T}(F)$ with the conditional expected sum over vertices of the completed ball $B_{\mathcal U}(0,q_n)$, with weight $b_n$ per vertex. Each vertex is recorded together with the first point $w$ where its path meets $\mathcal S$ and the full path from $w$ to the vertex. The additional atoms account for vertices already in $\mathcal S$. Repeated volume samples are included, and each new branch enlarges $\mathcal S$. Thus $\overline Q_{n,\mathcal S,T}$ describes, given $\mathcal S$, where the next volume-sampled vertex of the ball is attached to $\mathcal S$ and along which path.

\paragraph{Walk-step indices, LERW length, and cut times.}

Let $S$ be simple random walk from $0$, and give $\eta=\operatorname{LE}(S[0,\infty))$ its unconditioned law $\mathsf P$. Random-walk time counts steps of $S$, and the LERW length parameter counts edges of $\eta$. The following indices relate the two:
\[
\sigma_0=0,\qquad \sigma_{j+1}=1+\max\{k:S_k=S_{\sigma_j}\},\qquad \rho_m=\max\{j:\sigma_j\le m\},\qquad \eta_j=S_{\sigma_j}.
\]
Here $\sigma_j$ is the walk-step index selected for the $j$th vertex $\eta_j$ by the infinite loop-erasure construction. It need not be the last visit to $\eta_j$: that visit occurs at $\lambda_j=\sigma_{j+1}-1$. Thus $\rho_m$ counts the positive indices $\sigma_1,\sigma_2,\ldots$ that do not exceed $m$, or equivalently the number of edges in the prefix $\eta[0,\rho_m]$, which is contained in $\operatorname{LE}(S[0,m])$. It need not equal the length of this finite loop-erasure: the latter can have additional edges that are erased later. The cut-time estimates below compare these two lengths. The shift $\vartheta\eta=(\eta_{j+1}-\eta_1)_{j\ge0}$ deletes the first edge and translates the new starting point to $0$. The scales and the probability normalization of capacity are those of \eqref{min:scales} and \eqref{min:capacity-convention}. The capacity convention of \cite{HutchcroftSousi2023} differs by a factor of eight, which changes only the constants in the bounds below.

A time $v$ is a global cut time if $S[0,v]\cap S[v+1,\infty)=\varnothing$. By \cite[Lemma~7.7.4]{Lawler1991}, in the form recorded in \cite[Lemma~3.5]{HutchcroftSousi2023},
\begin{equation}
\mathbb P\{\text{no global cut time in }[a,b]\} \le C\frac{\log\log(b+3)}{\log(b+2)}, \qquad 3\le a<b,\quad b-a\ge b/(\log b)^6. \label{bridge:cutbound}
\end{equation}
At a global cut, erasing the rest of the walk cannot change the preceding loop-erasure. We use cuts near block endpoints to compare finite loop-erasures, the infinite loop-erasure, and loop-erasures of independent blocks of walk increments.

We also use the Gaussian maximal inequality \cite[Proposition~2.1.2, second assertion, equation~(2.7)]{LawlerLimic2010}:
\begin{equation}
\Pp_x\!\left\{\max_{0\le\ell\le m}|W_\ell-x|\ge a\right\} \le C\exp(-a^2/(Cm)),\qquad m\ge1,\quad a\ge0. \label{bridge:gaussian-maximal}
\end{equation}
For any integer $r_0\ge0$ and any event $H\in\mathcal F_{r_0}^W$, independence of the subsequent increments gives the weighted form
\[
\Pp_x\!\left( H,\max_{0\le\ell\le m}|W_{r_0+\ell}-W_{r_0}|\ge a\right) \le C\Pp_x(H)e^{-a^2/(Cm)}.
\]
When $H$ is an early avoidance event, this form keeps its probability as a factor in the displacement bound. We need that factor when summing errors over attachment points.

\begin{lemma}[Intersection of independent walk segments]
\label{prelim:independent-intersection}
Let $S,W$ be independent simple random walks from $0$. There is a constant $C$ such that, for all integers $a,b,m\ge0$ with $m+a\ge1$,
\begin{equation}\label{prelim:independent-intersection-bound}
\Pp\bigl(W[m,\infty)\cap S[a,a+b]\ne\varnothing\bigr) \le\frac{C(b+1)}{(m+a)\log(b+2)}.
\end{equation}
\end{lemma}
\begin{proof}
For a deterministic finite set $A$, last-exit decomposition \cite[Proposition~4.6.4]{LawlerLimic2010} gives
\[
\Pp\bigl(W[m,\infty)\cap A\ne\varnothing\bigr) =\sum_{x\in A}e_A(x)\sum_{t\ge m}p_t(0,x).
\]
Put $R_b=S[a,a+b]-S_a$. The increment range $R_b$ is independent of $S_a$ and has the law of $S[0,b]$. Conditional on $S_a$ and $R_b$, apply the preceding identity to $A=S_a+R_b$. Averaging first over $S_a$, symmetry, convolution and the heat-kernel bound \cite[Theorem~2.1.3]{LawlerLimic2010} give, for each $x\in R_b$,
\[
\Ee\!\left[\sum_{t\ge m}p_t(0,S_a+x)\,\middle|\,R_b\right] =\sum_{t\ge m}p_{t+a}(0,x)\le\frac{C}{m+a}.
\]
Translation invariance of capacity therefore gives an upper bound $C\Ee\operatorname{cap}(R_b)/(m+a)$. Finally, \cite[Corollary~1.4]{AsselahSchapiraSousi2018} gives $\Ee\operatorname{cap}(R_b)\le C(b+1)/\log(b+2)$, with the finitely many small values of $b$ covered by increasing $C$.
\end{proof}

\paragraph{The relation between walk steps and LERW length.}

With the testing horizon measured in random-walk steps, \cite[Theorem~1.1]{LawlerSunWu2019} gives the limit
\begin{equation}\label{prelim:raw-escape-normalization}
X_m=(\log m)^{1/3}\Pp_0^W(W[1,m]\cap\eta=\varnothing\mid\eta) \longrightarrow X_\infty, \qquad \alpha=\mathsf E_{\mathsf P}X_\infty\in(0,\infty),
\end{equation}
almost surely and in every finite $L^p(\mathsf P)$. Here $W$ is independent of the generating walk $S$, and $m$ counts steps of $W$. The source uses horizon $h^2$ and prefactor $(\log h)^{1/3}$, so $X_\infty$ is $2^{1/3}$ times its limit. Monotonicity extends the conversion from square horizons to all integer $m$. The proof below compares independent blocks and uses cut times near their endpoints, as in \cite[proof of Theorem~7.7.5]{Lawler1991} and \cite[proof of Proposition~3.4]{HutchcroftSousi2023}.

The heuristic is simple. Step $k$ is the last visit to a vertex of $\eta$ exactly when the future $S[k+1,\infty)$ avoids $\operatorname{LE}(S[0,k])$. After the first $k$ steps are reversed, this is the event that an independent walk avoids a LERW, which by \eqref{prelim:raw-escape-normalization} has probability about $\alpha(\log k)^{-1/3}$. Summing over $k<N$ gives $\rho_N\approx\alpha N/L_N$, and inverting gives $\sigma_m\approx mL_m/\alpha$. The lemma makes this quantitative.

\begin{lemma}[Walk steps needed for a given LERW length]\label{prelim:calibrated-clock}
For every $\varepsilon>0$ there is a finite $C_\varepsilon$ such that
\begin{equation}\label{bridge:ordinaryclock}
\Pp\!\left(\left|\frac{\alpha\sigma_m}{mL_m}-1\right|>\varepsilon\right) +\Pp\!\left(\left|\frac{L_m\rho_m}{\alpha m}-1\right|>\varepsilon\right) \le C_\varepsilon\frac{\log\log(m+3)}{L_m^2},\qquad m\ge3.
\end{equation}
\end{lemma}
\begin{proof}
We first compute the mean of $\rho_N$ and compare $\rho_N$ with a sum of independent loop-erased block lengths. The relation $\{\sigma_m\le N\}=\{\rho_N\ge m\}$ then gives the bound for $\sigma_m$.

\emph{Step 1: comparison of escape probabilities.} Let $S,W$ be independent simple random walks from $0$ and put
\[
a_m=\Pp(W[1,m]\cap\eta=\varnothing),\qquad q_m=\Pp(W[1,\infty)\cap\operatorname{LE}(S[0,m])=\varnothing).
\]
For $m\ge3$, the two probabilities differ by at most
\begin{equation}\label{prelim:raw-finite-escape-comparison}
|q_m-a_m|\le C\frac{\log\log(m+3)}{\log(m+2)}.
\end{equation}
Lemma~\ref{prelim:independent-intersection}, with the names of the walks interchanged for the first term, gives
\[
\Pp(W[0,m]\cap S[\lfloor m/2\rfloor,\infty)\ne\varnothing) +\Pp(W[m+1,\infty)\cap S[0,m]\ne\varnothing) \le C/\log(m+2).
\]
On the event of a global cut in $[\lceil m/2\rceil,m]$, finite and infinite loop-erasure share the prefix through that cut, and their symmetric difference is contained in $S[\lfloor m/2\rfloor,\infty)$. Apply the preceding bounds and \eqref{bridge:cutbound} to compare the two avoidance events, proving \eqref{prelim:raw-finite-escape-comparison}.

\emph{Step 2: the mean asymptotic for $\rho_N$.} For each walk-step index $k$, let $Y_k$ indicate that $k$ is the last visit to one of the infinite LERW vertices. Then
\[
Y_k:=\one_{\{k=\lambda_j\text{ for some }j\ge0\}} =\one_{\{\operatorname{LE}(S[0,k])\cap S[k+1,\infty)=\varnothing\}}.
\]
If the future avoids the finite loop-erasure, that loop-erasure is a prefix of $\eta$ and its endpoint is last visited at $k$. Conversely, the last-visit construction gives this avoidance whenever $k=\lambda_j$. Reverse the first $k$ walk steps and translate $S_k$ to $0$. Fixed-time loop-erasure reversibility \cite[Theorem~4.6]{HutchcroftSousi2023} and independence of future increments then give $\Ee Y_k=q_k$. Equations~\eqref{prelim:raw-escape-normalization} and \eqref{prelim:raw-finite-escape-comparison} imply $L_mq_m\to\alpha$. Since $\sigma_{j+1}\le N$ exactly when $\lambda_j<N$, summation of the slowly varying sequence $(q_m)$ gives
\begin{equation}\label{prelim:retained-clock-mean}
\rho_N=\sum_{k=0}^{N-1}Y_k, \qquad \bar\rho_N:=\Ee\rho_N=\sum_{k=0}^{N-1}q_k \sim\frac{\alpha N}{L_N}.
\end{equation}

\emph{Step 3: comparison with independent blocks.} For large $N$, divide $[0,N]$ into $M_N=\lceil(\log N)^2\rceil$ deterministic blocks $[t_{N,i-1},t_{N,i}]$ of lengths $B_{N,i}\asymp N/(\log N)^2$. Put
\[
\xi_i=|\operatorname{LE}(S[t_{N,i-1},t_{N,i}])|, \qquad h_N=\left\lceil4N/(\log N)^6\right\rceil.
\]
The variables $\xi_i$ are independent, since their lengths depend only on separate increment blocks, and $0\le\xi_i\le B_{N,i}$. For $t_{N,i-1}\le k<t_{N,i}$, the finite-block indicators are
\[
Y_k^{(i)} =\one_{\{\operatorname{LE}(S[t_{N,i-1},k]) \cap S[k+1,t_{N,i}]=\varnothing\}}, \qquad \xi_i=\sum_{k=t_{N,i-1}}^{t_{N,i}-1}Y_k^{(i)}.
\]
Call block $i$ good when there is a global cut in each of
\[
[t_{N,i-1}+\lceil h_N/2\rceil,t_{N,i-1}+h_N], \qquad[t_{N,i}-h_N,t_{N,i}].
\]
For large $N$, both windows lie inside the block, and their widths are at least $2N/\log^6N+O(1)$ and their endpoint logarithms are $\log N+O(\log\log N)$. Thus \eqref{bridge:cutbound} bounds the bad-block probability by $\Pp(\mathcal B_i)\le C\log\log N/\log N$. On a good block, let $v_1<v_2$ be the two cuts. No intersection crosses either cut. Consequently, for $v_1\le k\le v_2$ and $k<t_{N,i}$, both $Y_k$ and $Y_k^{(i)}$ equal
\[
\one_{\{\operatorname{LE}(S[v_1,k]) \cap S[k+1,v_2]=\varnothing\}}.
\]
They can disagree at no more than $2h_N+O(1)$ walk-step indices. Consequently
\[
\left|\rho_N-\sum_i\xi_i\right| \le CM_N(h_N+1)+C\sum_iB_{N,i}\one_{\mathcal B_i},
\]
Taking expectations, and using independence for the variance, gives
\begin{equation}\label{prelim:clock-block-errors}
\Ee\left|\rho_N-\sum_i\xi_i\right| \le CN\frac{\log\log N}{\log N}, \qquad \operatorname{Var}\left(\sum_i\xi_i\right) \le\sum_iB_{N,i}^2\le\frac{CN^2}{(\log N)^2}.
\end{equation}

\emph{Step 4: concentration and inversion.} The expectation bound also controls $|\bar\rho_N-\Ee\sum_i\xi_i|$. Markov's inequality for the comparison error and Chebyshev's inequality for the block sum yield
\[
\Pp\{|\rho_N-\bar\rho_N|>\varepsilon N/L_N\} \le C_\varepsilon\frac{\log\log(N+3)}{L_N^2}.
\]
The mean asymptotic \eqref{prelim:retained-clock-mean} replaces $\bar\rho_N$ by $\alpha N/L_N$ for all sufficiently large $N$, and increasing $C_\varepsilon$ covers the remaining values. Finally, $\sigma_m\le N$ if and only if $\rho_N\ge m$. Apply the resulting bound at $N_m^\pm=\lfloor(1\pm\delta)mL_m/\alpha\rfloor$, with fixed $0<\delta<1$, and use $L_{N_m^\pm}/L_m\to1$ to obtain the assertion for $\sigma_m$. This proves \eqref{bridge:ordinaryclock}.
\end{proof}

\paragraph{The Brownian limit for a finite loop-erased path.}

The occupation formula requires both a branch path and its rescaled number of edges. The following limit records both. Its normalization follows from Lemma~\ref{prelim:calibrated-clock} and the simple random-walk invariance principle.

\begin{corollary}[Brownian limit with the LERW lifetime recorded]
\label{bridge:finiteFCLT}
Let $S$ be simple random walk from $0$. For $t>0$, divide the spatial coordinates of $\operatorname{LE}(S[0,\lfloor tN_n\rfloor])$ by $R_n$ and give each edge duration $1/n$, recording the resulting lifetime. Its law converges to that of $B^{[t]}=(B_{s\wedge t})_{s\ge0}$, with lifetime $t$ and $\Ee B_tB_t^{\mathsf T}=tI_4/4$, uniformly for $t$ in compact subsets of $(0,\infty)$.
\end{corollary}
\begin{proof}
For every fixed $M<\infty$, \eqref{bridge:ordinaryclock}, monotonicity and a finite grid give
\[
\sup_{0\le s\le M} \left|\sigma_{\lfloor ns\rfloor}/N_n-s\right| \longrightarrow0\quad\hbox{in probability}.
\]
Donsker's theorem for $R_n^{-1}S_{\lfloor N_n\cdot\rfloor}$ and this time change give $R_n^{-1}\eta_{\lfloor n\cdot\rfloor}\Rightarrow B$, locally uniformly, with covariance $I_4/4$.  Linear interpolation changes the LERW path by at most $1/R_n$. Fix $0<\varepsilon<M<\infty$, put $k_{n,t}=\lfloor tN_n\rfloor$ and $h_{n,t}=\lceil k_{n,t}/(\log k_{n,t})^6\rceil$. By \eqref{bridge:cutbound}, uniformly for $t\in[\varepsilon,M]$, a global cut $v\in[k_{n,t}-h_{n,t},k_{n,t}]$ exists with probability tending to one, and then
\[
\operatorname{LE}(S[0,k_{n,t}]) =\eta[0,\rho_v]\oplus\operatorname{LE}(S[v,k_{n,t}]).
\]
The bound on $\rho_v$ gives $\rho_v/n-t\to0$ in probability. The final segment has at most $h_{n,t}=O_{\varepsilon,M}(nL_n^{-17})=o(n)$ edges, and its diameter divided by $R_n$ tends to zero by \eqref{bridge:gaussian-maximal} on the deterministic interval $[k_{n,t}-h_{n,t},k_{n,t}]$. This proves the finite stopped-path and lifetime convergence. The same argument applies to every sequence $t_n\in[\varepsilon,M]$ after passing to a subsequence with $t_n\to t$. Then, the stated uniformity follows from the continuity of $t\mapsto\Law(B^{[t]})$.
\end{proof}

\paragraph{A sharper upper tail for the walk-step indices.}

For the moment bounds below, we improve the upper tail for $\sigma_m$ in \eqref{bridge:ordinaryclock} to a bound proportional to $\log\log(m+3)/\log(m+2)$. We compare two independent blocks. When cuts occur near their endpoints, $\rho_N<m$ forces both blocks to have short loop-erasures. This is the argument of \cite[proof of Proposition~3.4]{HutchcroftSousi2023}.

\begin{lemma}[Upper tail for the walk-step index]\label{bridge:inverse-clock}
There are fixed finite constants $C_{\mathrm{time}},C$ such that, for all $m\ge3$,
\begin{equation}
\Pp(\sigma_m>C_{\mathrm{time}}mL_m) \le C\frac{\log\log(m+3)}{\log(m+2)}. \label{bridge:eqclock}
\end{equation}
\end{lemma}
\begin{proof}
For sufficiently large $m$, put $B_m=\lceil8mL_m/\alpha\rceil$ and $w_m=\lfloor m/100\rfloor$. Let $\mathcal C$ be the event that there is a global cut in each of
\[
[B_m-2w_m,B_m-w_m],\qquad [B_m+w_m,B_m+2w_m],\qquad [2B_m-2w_m,2B_m-w_m].
\]
By \eqref{bridge:cutbound}, $\Pp(\mathcal C^c)\le C\log\log m/\log m$. On $\mathcal C\cap\{\rho_{2B_m}<m\}$, the loop-erasure of $S[0,B_m]$ consists of a prefix of $\eta$ and a final segment of at most $2w_m+1$ steps. The loop-erasure of $S[B_m,2B_m]$ consists of a segment of $\eta$ between its two cuts and two end segments containing at most $4w_m+2$ steps in total. The segments belonging to $\eta[0,\rho_{2B_m}]$ have fewer than $m$ edges. Therefore
\[
\{\rho_{2B_m}<m\}\cap\mathcal C \subseteq\bigcap_{i=0}^1 \{|\operatorname{LE}(S[iB_m,(i+1)B_m])|\le m+4w_m+2\}.
\]
The two block events are independent. Extend each block independently to an infinite walk. For a block of $B_m$ steps, its loop-erasure contains the first $\rho_{B_m}$ edges of its infinite loop-erasure. Thus \eqref{bridge:ordinaryclock} and $\alpha B_m/L_{B_m}\sim8m$ bound each probability by $C\log\log m/(\log m)^{2/3}$. Consequently
\[
\Pp(\rho_{2B_m}<m) \le C\frac{\log\log m}{\log m} +C\frac{(\log\log m)^2}{(\log m)^{4/3}} \le C\frac{\log\log m}{\log m}.
\]
Since $\{\sigma_m>2B_m\}=\{\rho_{2B_m}<m\}$ and $2B_m\le C_{\mathrm{time}}mL_m$ for a fixed $C_{\mathrm{time}}$, this proves the claim; increasing $C$ covers small $m$.
\end{proof}

\paragraph{Local escape probabilities and uniform integrability.}

For later occupation estimates we need moment bounds for sums of local escape probabilities. We obtain them from the capacity of random-walk paths and the upper tail for $\sigma_m$.

Fix $0<T<\infty$ and $q_n=\lfloor nT\rfloor$. The three cutoffs are
\begin{equation}\label{not:cutoffs}
c_n=1\vee\lfloor nL_n^{-40}\rfloor,\qquad r_n=1\vee\lfloor nL_n^{-12}\rfloor,\qquad u_n=1\vee\lfloor nL_n^{-6}\rfloor.
\end{equation}
Here $c_n$ is a window radius in LERW edges, while $r_n,u_n$ count random-walk steps. The choice $c_n\ll r_n\ll u_n$ lets us first replace avoidance for all time by avoidance through step $r_n$, then separate those first $r_n$ steps from the walk after step $u_n$. All three are much smaller than $n$, so on the scale of the theorem they are invisible. The particular powers of $L_n$ are chosen for convenience: we only use that $c_n$, $r_n$, $u_n$ and $n$ are separated by sufficiently large powers of $L_n$. For $0\le j\le q_n$ define
\[
\mathcal C_j^\eta=\eta[(j-c_n)_+,j+c_n],\quad e_j=\Pp_{\eta_j}(\tau_{\mathcal C_j^\eta}^+=\infty),\quad p_j=\Pp_{\eta_j}(W[1,r_n]\cap \mathcal C_j^\eta=\varnothing),
\]
where all probabilities in these definitions are conditional on $\eta$. Put
\[
Z_{n,T}=\frac{L_n^2}{n}\sum_{j=0}^{q_n}e_j, \qquad E_{n,T}=\frac{L_n^2}{n}\sum_{j=0}^{q_n}p_j.
\]
The normalization $L_n^2/n$ is the capacity scale of \eqref{min:kappa-normalization}, so both sums are of order one. Capacities bound $Z_{n,T}$. The events defining $E_{n,T}$ use only the first $r_n$ test-walk steps, as needed in the later stopped-path comparison.
\begin{lemma}[Moments of the local escape-probability sums]\label{bridge:rows}
For every $1\le p<3/2$,
\begin{equation}
\sup_n\Ee[Z_{n,T}^p+E_{n,T}^p]<\infty. \label{bridge:eqrowLp}
\end{equation}
\end{lemma}
\begin{proof}
We bound $Z_{n,T}$ by capacities of random-walk paths, then compare $E_{n,T}$ directly with $Z_{n,T}$.

Choose a deterministic $N_n^*\le C_TnL_n$ so that Lemma~\ref{bridge:inverse-clock} implies
\[
\Pp(\sigma_{q_n}>N_n^*)\le C_T\log\log n/L_n^3.
\]
Partition the walk-step interval $[0,N_n^*]$ into consecutive blocks of $c_n$ steps, extending the last one to a full block. Let $V_b$ contain those LERW vertices $\eta_j$ whose associated indices $\sigma_j$ lie in block $b$. Assign a vertex at a shared endpoint to just one block. The LERW indices in $V_b$ are consecutive and span at most $c_n$ edges, because the $\sigma_j$ increase by at least one per edge. Thus $V_b\subset \mathcal C_j^\eta$ for every $j$ assigned to this block. On $\{\sigma_{q_n}\le N_n^*\}$, monotonicity of escape probabilities and of capacity gives
\[
\sum_{j=0}^{q_n}e_j \le\sum_b\sum_{x\in V_b}e_{V_b}(x) =\sum_b\operatorname{cap}(V_b) \le\sum_b\operatorname{cap}(S[bc_n,(b+1)c_n]).
\]
Capacity is translation invariant and the blocks use independent increments, so these capacities are independent and identically distributed. Set
\[
Y_n=\frac{L_n^2}{n} \sum_{b<\lceil N_n^*/c_n\rceil}\operatorname{cap}(S[bc_n,(b+1)c_n]).
\]
The random-walk capacity estimate of \cite[Corollary~1.4]{AsselahSchapiraSousi2018}, recorded as \cite[Lemma~2.5]{HutchcroftSousi2023}, gives $\Ee\operatorname{cap}(S[0,c_n])\le Cc_n/\log(c_n+2)$.  Since $\operatorname{cap}(S[0,c_n])\le c_n+1$,
\[
\Ee Y_n\le C_T, \qquad \operatorname{Var}(Y_n) \le C_T\frac{L_n^5c_n}{n}=O_T(L_n^{-35}).
\]
Thus $Y_n$ is bounded in $L^2$.  Always $Z_{n,T}\le C_TL_n^2$, so
\[
\Ee[Z_{n,T}^p;\sigma_{q_n}>N_n^*] \le C_{T,p} L_n^{2p-3}\log\log n\longrightarrow0 \qquad(p<3/2).
\]
On the complementary event $Z_{n,T}\le Y_n$, so the $L^2$ bound for $Y_n$ proves the required $L^p$ bound for $Z_{n,T}$.

For avoidance through step $r_n$, the heat-kernel bound $p_k(x,y)\le Ck^{-2}$ gives
\[
0\le p_j-e_j \le\Pp_{\eta_j}(W[r_n,\infty)\cap \mathcal C_j^\eta\ne\varnothing) \le C(c_n+1)/r_n.
\]
Consequently
\begin{equation}
0\le E_{n,T}-Z_{n,T}\le C_TL_n^2c_n/r_n+o(1)=O_T(L_n^{-26})+o(1). \label{bridge:earlyrowerror}
\end{equation}
The deterministic comparison \eqref{bridge:earlyrowerror} transfers this bound to $E_{n,T}$.
\end{proof}

\paragraph{The contribution from walks taking many steps.}

A walk can take far more than $N_n$ steps and still produce a short loop-erasure, by spending most of its time in loops that are erased later. The occupation identity counts every step, so we need to show that such walks contribute little. The next bounds combine Lemma~\ref{bridge:inverse-clock}, fixed-time reversal \cite[Theorem~4.6]{HutchcroftSousi2023}, and the conditional duration estimate in \cite[Lemma~2.6 and equation~(9)]{HalberstamHutchcroft2024}. That estimate holds for every absorbing set.

The quantity $H_{A,a}(J,N)$ counts, in expectation, the step indices $k>N$ before the walk returns to $A$ for which the loop-erasure has between $1$ and $J$ edges. The bound is uniform in $A$, so it still applies after further branches have been exposed.

\begin{corollary}[Occupation bounds by walk duration and LERW length]
\label{closure:extreme}
For a set $A\subseteq\mathbb Z^4$ containing $a$, kill the testing walk on its first return to $A$ and put
\[
H_{A,a}(J,N)=\mathbb E_a\sum_{k>N} \mathbf1_{\{k<\tau_A^+,\;1\le |\operatorname{LE}(W[0,k])|\le J\}}.
\]
There are fixed constants $C,C_{\mathrm{dur}}$ such that, for every integer $J\ge1$,
\begin{equation}\label{bridge:universalH}
H_{A,a}(J,0)\le CJ(\log(J+2))^{1/3},\qquad H_{A,a}(J,N)\le\frac{CJ}{N}H_{A,a}(J,0) \quad\bigl(N\ge C_{\mathrm{dur}}J\log(J+2)\bigr).
\end{equation}
The second inequality also holds after weighting each loop-erased path by any nonnegative function of that path.
\end{corollary}
\begin{proof}
We first condition on the loop-erased path and bound the number of walk steps. Summing over paths gives the tail bound. We then use the upper tail for $\sigma_m$ to bound total occupation.

For each admissible simple path $\lambda=(a,x_1,\ldots,x_\ell)$, $1\le\ell\le J$, set $m_\lambda(k)=\Pp_a(k<\tau_A^+,\operatorname{LE}(W[0,k])=\lambda)$. Fixed-time reversibility \cite[Theorem~4.6]{HutchcroftSousi2023} identifies this measure on numbers of walk steps with the hitting-time measure of a reversed walk from $x_\ell$ that first hits the full set $A$ at $a$ and has loop-erasure $\lambda^\leftarrow$. The conditional-duration estimate holds for infinite $A$ as well \cite[Lemma~2.6]{HalberstamHutchcroft2024}, so the estimate is uniform in the entire exposed tree. By the same lemma, its conditional hitting time has a centering $T_{A,\varnothing} (\lambda^\leftarrow)$ and tail probability at most $C\ell/v$ at distance $v\ge\ell$ from that centering. The bound in \cite[equation~(9)]{HalberstamHutchcroft2024}, with each escape probability at most one, gives
\[
T_{A,\varnothing}(\lambda^\leftarrow) \le C\ell\sum_{s=1}^{\ell}\frac1s \le C\ell\log(\ell+1).
\]
 Hence, for $N\ge C_{\mathrm{dur}}J\log(J+2)$ with $C_{\mathrm{dur}}$ sufficiently large,
\[
\sum_{k>N}m_\lambda(k) \le\frac{C\ell}{N}\sum_{k\ge0}m_\lambda(k) \le\frac{CJ}{N}\sum_{k\ge0}m_\lambda(k).
\]
Multiplying by an arbitrary nonnegative path weight and summing proves the weighted tail bound. To bound the total occupation, it suffices by monotonicity to take $A=\{a\}$ and write $H(J,N)=H_{\{a\},a}(J,N)$. This quantity is finite: a loop-erased path of at most $J$ edges ends in the finite lattice ball $B_{\mathbb Z^4}(a,J)$, so
\[
H(J,0)\le\sum_{x\in B_{\mathbb Z^4}(a,J)}G(a,x)<\infty.
\]
Choose a fixed $C_{\mathrm{blk}}\ge C_{\mathrm{dur}}$ and put $M_J=\lceil C_{\mathrm{blk}}J\log(J+2)\rceil$, with $C_{\mathrm{blk}}$ large enough that the established tail bound gives $H(J,M_J)\le H(J,0)/2$. Consequently
\[
H(J,0)\le2\sum_{k=1}^{M_J}\Pp(J_k\le J).
\]
Extend the testing walk $W$ independently to infinity. In this comparison, let $\eta$ be its infinite loop-erasure and let $\sigma,\rho$ be the indices defined from this extended walk. The loop-erasure through step $k$ contains $\eta[0,\rho_k]$, so $J_k\ge\rho_k$ and $\{J_k\le J\}\subseteq\{\sigma_{J+1}>k\}$. For $k>C_{\mathrm{time}}(J+1)L_{J+1}$, Lemma~\ref{bridge:inverse-clock} bounds this probability by $C\log\log(J+3)/\log(J+2)$. Splitting the last sum at this threshold, and enlarging constants for small $J$, yields
\[
H(J,0)\le CJ L_J+CJ\log\log(J+3)\le CJ L_J.
\]
This proves the first bound for every absorbing set $A$.
\end{proof}

\subsection{The spatial Brownian Kesten tree}\label{prelim:continuum}

We fix the excursion normalization, construct the limiting tree, and identify its ICSBM occupation measure and ball-volume law. The spinal decomposition describes the tree as seen from a sampled vertex and yields a recursion for successive volume samples. We then state the moment criterion for identifying the limit. Throughout this subsection, $\kappa>0$ is the capacity constant in Section~\ref{min:setup} and $\Sigma=I_4/4$.

\paragraph{Brownian excursions and the Kesten construction.}
Let $\mathbb N_x^{\kappa,\Sigma}$ be the spatial Brownian-tree excursion measure with branching mechanism $\psi(q)=\kappa q^2/2$ and spatial motion a Brownian motion of covariance $\Sigma$. Its level measures $\ell^t$ and volume measure $m$ are normalized by
\begin{equation}\label{min:snake-normalization}
m(\cdot)=\int_0^\infty\ell^t(\cdot)\,dt, \qquad \mathbb N_x^{\kappa,\Sigma} [1-e^{-\lambda\|\ell^t\|}] =\frac{\lambda}{1+\kappa t\lambda/2}.
\end{equation}
Here $\ell^t$ is carried by vertices at intrinsic height $t$. These are the critical quadratic Ray-Knight normalizations \cite[Theorem~1.4.1 and Corollary~1.4.2]{DuquesneLeGall2002}. Write $H$ for the maximum distance from the excursion root. In particular, $\mathbb N_x[\|\ell^t\|]=1$ and $\mathbb N_x(H>t)=2/(\kappa t)$ for $t>0$. We omit the superscripts on $\mathbb N_x$ when the parameters are fixed.

The underlying measured Brownian Kesten tree is the critical quadratic specialization of the size-biased immigration L\'evy tree \cite[Sections~1.2--1.3 and~3.3]{Duquesne2009}. Its spatial marking is the Brownian-snake marking of \cite[Section~6]{DuquesneLeGall2005}. In our normalization, $\mathcal K_{\kappa,\Sigma}$ consists of an infinite spine $[0,\infty)$, rooted at $0$ and labelled by a $\Sigma$-Brownian motion $\xi$, together with conditionally independent spatial excursions grafted according to the Poisson intensity
\begin{equation}\label{min:kesten-ppp}
\kappa\,ds\,\mathbb N_{\xi_s}(d\mathcal E).
\end{equation}
Its metric is the grafted tree metric, and its volume is the sum of the excursion volumes, with zero volume on the spine. On each bounded spine interval there are finitely many grafts of height exceeding any fixed positive number. Campbell's formula  (see, e.g., \cite[Proposition~2.7]{LastPenrose2018}) and \eqref{min:snake-normalization} also give
\begin{equation}\label{min:kesten-first-moment}
\Ee[m_{\mathcal K}(B(o,T))] =\kappa\int_0^T\int_0^{T-s}dt\,ds =\frac{\kappa T^2}{2}.
\end{equation}
Consequently its volume is locally finite almost surely. The level-measure disintegration in \eqref{min:snake-normalization} makes its radial volume measure absolutely continuous with respect to $dt$. Thus, every intrinsic sphere has zero volume almost surely. This is the boundary condition needed when restricting the limiting measured tree to intrinsic balls. We record the following connection between Brownian Kesten tree and ICSBM.

\begin{proposition}[ICSBM occupation measure and ball volume]\label{prelim:icsbm-volume}
Write $\mathcal K=\mathcal K_{\kappa,\Sigma}$, let $\phi$ be its spatial label map, and put $h(u)=d_{\mathcal K}(o,u)$. For the grafts $(s_i,\mathcal E_i)$ in \eqref{min:kesten-ppp}, let $\phi_{\mathcal E_i}$ be their spatial label maps and define
\[
\mathcal X_t=\sum_{s_i<t}(\phi_{\mathcal E_i})_*\ell_{\mathcal E_i}^{\,t-s_i}\quad(t>0),\qquad \mathcal X_0=0.
\]
Then $\mathcal X$ is ICSBM issued from the origin, with branching mechanism $\psi(q)=\kappa q^2/2$ and spatial covariance $\Sigma$, and
\begin{equation}\label{min:kesten-occupation}
(h,\phi)_*m_{\mathcal K}=dt\,\mathcal X_t(dx).
\end{equation}
For every $T>0$, $V_T=m_{\mathcal K}(B(o,T))=\int_0^T\|\mathcal X_t\|\,dt$ has Laplace transform
\begin{equation}\label{min:ball-volume-laplace}
\mathbb E e^{-\lambda V_T}=\cosh^{-2}\!\left(T\sqrt{\frac{\kappa\lambda}{2}}\right)=\cosh^{-2}\!\left(\pi T\sqrt{\frac{\lambda}{24\alpha}}\right),\qquad\lambda\ge0.
\end{equation}
\end{proposition}
\begin{proof}
The identification of $\mathcal X$ is the immortal-particle representation \cite{Evans1993}, in its canonical form \cite[Section~3.7]{VanDerHofstad2006}. Disintegrating each excursion volume by \eqref{min:snake-normalization} and recalling that the spine has zero volume gives \eqref{min:kesten-occupation}.

For $\lambda>0$, put $u_\lambda(t)=\mathbb N_0[1-\exp\{-\lambda\int_0^t\|\ell^s\|\,ds\}]$. The log-Laplace equation for the branching mechanism $\psi(q)=\kappa q^2/2$ gives $u_\lambda'=\lambda-\kappa u_\lambda^2/2$, with $u_\lambda(0)=0$, and hence $u_\lambda(t)=\sqrt{2\lambda/\kappa}\tanh(t\sqrt{\kappa\lambda/2})$. The Poisson construction \eqref{min:kesten-ppp} yields $\mathbb E e^{-\lambda V_T}=\exp\{-\kappa\int_0^T u_\lambda(t)\,dt\}$, proving \eqref{min:ball-volume-laplace}; the second equality uses \eqref{min:escape-constant-check}, and the case $\lambda=0$ follows by continuity.
\end{proof}

\paragraph{The spinal decomposition.}
For a vertex $u$ at height $t$, write $W_u(s)$ for the spatial label of its ancestor at height $s$, $0\le s\le t$. Under the level-$t$ Campbell probability $\mathbb N_x(d\mathcal E)\ell_{\mathcal E}^t(du)$, whose total mass is one, the ancestral labels form a $\Sigma$-Brownian path of lifetime $t$. Conditional on that path, the components off its ancestral segment, with planar order forgotten, form a Poisson collection with intensity $\kappa\,ds\,\mathbb N_{W_u(s)}$. If $\mathcal P_u$ records these grafts and $F,g\ge0$, the decomposition reads
\begin{align}\label{min:spinal-formula}
&\mathbb N_x\!\left[ \int\ell^t(du)F(W_u) \exp\{-\langle\mathcal P_u,g\rangle\}\right]\notag\\
&\quad=\Ee_x^\Sigma\!\left[ F(B_{[0,t]}) \exp\left\{-\kappa\int_0^t \mathbb N_{B_s}(1-e^{-g(s,\mathcal E)})\,ds\right\}\right].
\end{align}
Formula \eqref{min:spinal-formula} is the critical quadratic specialization of the ancestral decomposition \cite[Theorem~4.5]{DuquesneLeGall2005}, with the independent spatial marking of \cite[Section~6]{DuquesneLeGall2005}; see also \cite[Lemmas~4.2.4--4.2.5 and~4.6.1]{DuquesneLeGall2002} for its snake formulation. The two contour sides each have intensity $(\kappa/2)\,ds\,\mathbb N_{W_u(s)}$, giving the displayed total intensity after planar order is forgotten. Preservation of this Poisson grafting form allows the same branch insertion rule to be iterated under the volume-integrated sampling measures below.

\paragraph{Successive volume sampling.}
Let $\mathcal T$ consist of the infinite spine together with finitely many oriented branches added successively at points of the preceding tree. Write $\lambda_{\mathcal T}$ for intrinsic length measure, $h_{\mathcal T}(w)=d_{\mathcal T}(o,w)$, and $\phi$ for its labels. For a nonnegative branch test $F$, set
\begin{equation}\label{min:continuum-row}
Q_{\mathcal T,T}(F)=\kappa \int_{\mathcal T\cap B(o,T)}\lambda_{\mathcal T}(dw) \int_0^{T-h_{\mathcal T}(w)} \Ee_{\phi(w)}^\Sigma[F(w,B_{[0,t]})]\,dt.
\end{equation}
Attaching the sampled path at $w$ enlarges the exposed tree and records its endpoint as the next ordered sample. This is the continuum counterpart of \eqref{min:wilson-row}: length measure replaces the sum over vertices of the exposed tree, while the Brownian path of lifetime $t$ replaces the finite loop-erasure parametrized by its number of edges. The condition $h_{\mathcal T}(w)+t\le T$ is the same intrinsic-ball restriction. Informally, a new volume sample is reached by leaving the tree spanned so far at a point chosen according to $\kappa$ times length, and then following a Brownian path whose duration $t$ is chosen according to Lebesgue measure, subject to ending within distance $T$ of the root.

The following is an ordered, volume-integrated corollary of the classical spinal decomposition \eqref{min:spinal-formula} and the Mecke equation \cite[Theorem~4.1]{LastPenrose2018}. We record the short induction for the recorded exposed tree and ordered samples, together with its factorial bound. The total mass of the $k$-point measure is the expected $k$th power of ball volume, so the bound also controls the moments needed to identify the limiting measured ball.

\begin{corollary}[Successive volume sampling]\label{min:continuum-campbell}
Let $\mathcal H_k$ consist of the infinite spine and the reduced tree spanned by $k$ ordered volume samples from $B_T(\mathcal K_{\kappa,\Sigma})$, with the root, sample order, attachment positions and spatial paths retained. Define its Campbell measure $\mathfrak C_{k,T}$ by integrating against $m_{\mathcal K}^{\otimes k}$ and then taking expectation, without normalizing the volume measure.  Let $\mathfrak C_{0,T}$ be the law of the labelled infinite spine. Then $\mathfrak C_{k+1,T}$ is obtained from $\mathfrak C_{k,T}$ by \eqref{min:continuum-row}. In particular,
\begin{equation}\label{min:continuum-factorial}
\Ee[m_{\mathcal K}(B(o,T))^k] \le k!\left(\frac{\kappa T^2}{2}\right)^k, \qquad k\ge0.
\end{equation}
\end{corollary}

\begin{proof}
By the Mecke equation and \eqref{min:spinal-formula}, sampling a vertex from volume replaces its excursion by the ancestral Brownian segment to that vertex and an independent Poisson collection of side excursions. The unselected excursions retain their law. Consequently, after each insertion the unexplored grafts have intensity $\kappa\,\lambda_{\mathcal T}(dw)\mathbb N_{\phi(w)}(d\mathcal E)$ along the exposed tree $\mathcal T$. Disintegrating volume by \eqref{min:snake-normalization} gives exactly \eqref{min:continuum-row}, with the restriction $h_{\mathcal T}(w)+t\le T$. Induction yields the asserted recursion for all nonnegative tests of the exposed tree and ordered samples.  The volume is diffuse by \eqref{min:snake-normalization} and supported on leaves \cite[Theorem~4.6(i)]{DuquesneLeGall2005}, so exposed branches carry zero volume. Their endpoints have zero length measure, and no additional insertion terms arise.

To bound the total mass, set $v_T=\kappa T^2/2$, the expected volume in \eqref{min:kesten-first-moment}. After $k$ insertions, $\mathcal T$ is the union of the original spine and $k$ added branches, apart from their shared endpoints. For an added branch of length $\ell\le T-a$ starting at height $a$, its contribution to $\|Q_{\mathcal T,T}\|$ is
\[
\kappa\int_0^\ell(T-a-s)\,ds\le v_T.
\]
The spine contributes $v_T$, so $\|Q_{\mathcal T,T}\|\le(k+1)v_T$. The recursion and $\|\mathfrak C_{0,T}\|=1$ now give $\|\mathfrak C_{k,T}\|\le k!v_T^k$, proving \eqref{min:continuum-factorial}. The resulting exponential moment for parameters smaller than $v_T^{-1}$ supplies the moment-determinacy condition used below.
\end{proof}

\paragraph{Two facts about the local topology.}
All marked convergence statements use the space of marked measures on $X\times\mathbb R^4$ described in Section~\ref{intro:main-result}; labelled objects whose marks are given by a function carry its induced topology, and the candidate limits below are of this kind. Restriction to the ball of radius $T$ is continuous at spaces with $m\{x:d(o,x)=T\}=0$. The argument for restricting to balls in \cite[Lemma~2.6]{AthreyaLohrWinter2016} applies with spatial marks: restricting a rooted distance and label array imposes only $d(o,x_i)\le R$ and preserves its labels. In particular, finite rooted marked Gromov-weak convergence implies local convergence. A finite measure assigns positive mass to at most countably many intrinsic spheres, so ball restriction is continuous at almost every radius of the limiting space, and dominated convergence applies to \eqref{min:local-metric}.

\paragraph{The rooted marked method of moments.}
The following rooted form of \cite[Theorem~4.2]{FoutelRodier2024} identifies the law of the measured ball from all its expected rooted marked distance-matrix measures. These measures record distances from the root and between sampled vertices, together with their spatial labels. The factorial volume bound above supplies the moment-determinacy condition.

\begin{lemma}[Rooted marked metric-measure method of moments]
\label{prelim:marked-moments}
Let $\mathsf X_n$ and $\mathsf X$ be random rooted finite spatial metric-measure spaces as defined in Section~\ref{min:setup}, with root label $0$. Let $\nu_k$ be the rooted distance-matrix measure in \eqref{min:distance-matrix}. Suppose that, for every $k\ge1$,
\[
\Ee[\nu_k^{\mathsf X_n}] \Longrightarrow \Ee[\nu_k^{\mathsf X}]
\]
weakly as finite measures on the rooted distance and label array space. Write $M=\|m_{\mathsf X}\|$ and suppose that
\[
\sum_{r\ge1} \big(\Ee[(1+M)^r]\big)^{-1/(2r)}=\infty.
\]
Then $\mathsf X_n\Longrightarrow\mathsf X$ in rooted marked Gromov-weak topology. In particular, the mass condition holds whenever $\Ee[M^r]\le r!C^r$ for all $r\ge1$ and some $C<\infty$.
\end{lemma}

\begin{proof}
Augment the mark space to $\mathbb R^4\times\{0,1\}$ and the sampling measure to
\[
\mu^\bullet=(\mathrm{id},\phi,0)_*m+\delta_{(o,0,1)}.
\]
The extra atom records the root using a color absent from the original measure. Expanding an augmented moment by its color pattern gives a finite sum of continuous pushforwards of the assumed rooted moment measures; the all-color-$1$ term is deterministic. Thus all augmented moments converge, including the zeroth moment, which is one. Their limiting sampling mass is $1+M$, so \cite[Theorem~4.2]{FoutelRodier2024} gives weak convergence of the augmented marked spaces.

The rooted $k$-point measure is recovered from the augmented $(k+1)$-point measure by restricting to the color pattern $(1,0,\ldots,0)$, which is both open and closed, and projecting the corresponding coordinates. Because this color pattern is both open and closed, recovery of the root and rooted distance-matrix measures is continuous on the subspace of augmented spaces carrying a unit color-$1$ root atom. The resulting convergence is therefore rooted. Finally, the factorial moment bound and the binomial expansion give
\[
\Ee[(1+M)^r]\le r!(1+C)^r\le[r(1+C)]^r,
\]
so the inverse $2r$th roots are bounded below by $[r(1+C)]^{-1/2}$ and their sum diverges.
\end{proof}

The remaining task is to obtain this recursion from the discrete Wilson insertion measures. We begin with the capacity and escape laws that determine their limiting attachment intensity.

\section{Convergence of the local escape measure}
\label{min:capacity-section}

Recall the ordinary LERW law $\mathsf P$, the centered shift $\vartheta$ and the local escape probabilities from Section~\ref{prelim:lerw}. We prove convergence of the normalized local escape measure to $\kappa\,ds$ by comparing local-window escape sums with finite-path capacity. We estimate their difference in expectation under the stationary path law, then transfer the estimate to ordinary LERW.

Fix $0<T<\infty$, put $q_n=\lfloor nT\rfloor$, and use the cutoffs $c_n,r_n$ from \eqref{not:cutoffs} (recall that $c_n=1\vee\lfloor nL_n^{-40}\rfloor$ is a window radius in LERW edges and $r_n=1\vee\lfloor nL_n^{-12}\rfloor$ is a number of walk steps). Thus
\[
\mathcal C_j^\eta=\eta[(j-c_n)_+,j+c_n],\qquad p_j=\Pp_{\eta_j}(W[1,r_n]\cap \mathcal C_j^\eta=\varnothing\mid\eta).
\]
The following proposition identifies the $\kappa\,\mathrm ds$ attachment factor in the limiting insertion kernel, under the ordinary law $\mathsf P$. Section~\ref{sec:bridge} will add the branch path to this attachment measure.

\begin{proposition}[Convergence of the local escape measure]\label{bridge:row-convergence}
As finite measures on $[0,T]$,
\begin{equation}
\frac{L_n^2}{n}\sum_{j=0}^{q_n}p_j\delta_{j/n} \longrightarrow\kappa\,\one_{[0,T]}(s)\,\mathrm ds \quad\hbox{in probability}. \label{bridge:rowmeasure}
\end{equation}
The total masses converge in every $L^p$, $1\le p<3/2$.
\end{proposition}

\subsection{The stationary path law}

We need escape probabilities from the whole infinite path and from the path with its initial vertex removed. For the whole infinite path, use $X_m,X_\infty$ with the normalization by the number of testing-walk steps in \eqref{prelim:raw-escape-normalization}. Put $\eta^+=\eta[1,\infty)$, let $\mathrm{Es}_m,\mathrm{Es}_m^+$ be the conditional probabilities that $W[1,m]$ avoids $\eta$ and $\eta^+$, respectively, and set $\mathrm{Es}_0=1$. Conditional on $\eta$, define
\[
v_k=\Pp_0^W(W_k=0,\ W[1,k]\cap\eta^+=\varnothing\mid\eta),\qquad G_\eta=\sum_{k\ge0}v_k=G_{\mathbb Z^4\setminus\eta^+}(0,0).
\]
Here the Green function counts visits before hitting $\eta^+$. Write $X_m^+=(\log m)^{1/3}\mathrm{Es}_m^+$ and $X_\infty^+=G_\eta X_\infty$.

To average along the path, introduce the escape-biased law
\begin{equation}
\frac{d\mathsf Q}{d\mathsf P} =f:=\frac{X_\infty}{\alpha}, \label{min:density}
\end{equation}
which is the positive-half law of stationary two-sided LERW $\widehat\eta=(\widehat\eta_j)_{j\in\mathbb Z}$ \cite[Section~6]{LawlerSunWu2019}. We also denote the two-sided law by $\mathsf Q$; functionals of nonnegative indices use its positive-half marginal. In particular, $\mathsf Q$ is $\vartheta$-stationary. In all these normalized escape probabilities, the factor $\log m$ is indexed by the number of testing-walk steps.

The following lemma transfers convergence from $\mathsf Q$ to $\mathsf P$ when fixed shifts have uniformly vanishing effect.

\begin{lemma}[Transfer from the stationary law]\label{min:support}
There exists a finite integer-valued random variable $J_\ast$ under $\mathsf P$ such that each measure
\[
\mu_k(B)=\mathsf P(J_\ast=k,\ \vartheta^k\eta\in B)
\]
is absolutely continuous with respect to $\mathsf Q$. Consequently, if measurable real path functionals $F_n$ satisfy $F_n\to a\in\mathbb R$ in $\mathsf Q$-probability and $\sup_\gamma|F_n(\vartheta^k\gamma)-F_n(\gamma)|\to0$ for every fixed $k$, then $F_n\to a$ in $\mathsf P$-probability.
\end{lemma}
\begin{proof}
First note that, for every fixed $k$,
\begin{equation}
\mathsf P\circ \vartheta^{-k}\ll\mathsf P. \label{min:nonsingular}
\end{equation}
Indeed, conditional on a feasible prefix $\eta[0,k]=\gamma$, the domain Markov property \cite[Section~7.3]{Lawler1991} identifies the tail with the loop erasure of an SRW from $\gamma_k$ conditioned to avoid $\gamma[0,k-1]$ (the empty set if $k=0$). This conditioning event has positive probability.  Its walk law, and hence its loop-erased law, is absolutely continuous with respect to the unconditioned law.  Summing over the countably many prefixes proves \eqref{min:nonsingular}.

\emph{A zero-one law under $\mathsf P$.} We prove a zero-one law for events invariant under $\vartheta$, modulo $\mathsf P$-null sets. For deterministic $m$, let $\eta^{(m)}=\operatorname{LE}(S[m,\infty))-S_m$. By \eqref{bridge:cutbound}, the probabilities that there is no global cut in $[\lceil e^{r^2}\rceil,2\lceil e^{r^2}\rceil]$ are summable in $r$. Borel-Cantelli therefore gives arbitrarily late global cuts almost surely. At a cut time after $m$, the two loop erasures $\eta$ and $\eta^{(m)}$ have the same translated infinite tail.  Thus there are finite random indices $i,j$ with $\vartheta^i\eta=\vartheta^j\eta^{(m)}$. For an invariant event $E$, \eqref{min:nonsingular} lets us impose invariance at every deterministic shift outside one null set. Since $\eta^{(m)}$ also has law $\mathsf P$, the common tail gives $\one_E(\eta)=\one_E(\eta^{(m)})$ almost surely. The right-hand side depends only on SRW increments after $m$. Reverse-martingale convergence and Kolmogorov's zero-one law therefore give $\mathsf P(E)\in\{0,1\}$.

\emph{Transfer to ordinary LERW.} Set $A=\{f>0\}$. We show that a finite shift brings the path into $A$. Stationarity and $\mathsf Q(A)=1$ imply $\mathsf P(A\setminus \vartheta^{-1}A)=0$.  By \eqref{min:nonsingular}, the sets $\vartheta^{-k}A$ are increasing modulo $\mathsf P$-null sets.  Consequently $H=\bigcup_{k\ge0}\vartheta^{-k}A$ is invariant modulo null sets.  Since $\mathsf P(H)\ge\mathsf P(A)>0$, the preceding zero-one law gives $\mathsf P(H)=1$.  Take $J_\ast=\inf\{k:\vartheta^k\eta\in A\}$. If $\mathsf Q(B)=0$, then $\mathsf P(B\cap A)=0$, and \eqref{min:nonsingular} yields $\mu_k(B)\le\mathsf P(\vartheta^k\eta\in B\cap A)=0$.

Finally fix $M$ and $\varepsilon>0$. Uniform vanishing of the fixed-shift errors gives, for all large $n$,
\[
\mathsf P(|F_n-a|>\varepsilon) \le\mathsf P(J_\ast>M)+\sum_{k=0}^M\mu_k(|F_n-a|>\varepsilon/2).
\]
The sum tends to zero by $\mu_k\ll\mathsf Q$. The claim follows by first letting $n\to\infty$ and then letting $M\to\infty$.
\end{proof}

\subsection{Capacity per unit LERW length}

Each capacity increment is a product of the escape probabilities with and without the initial vertex \cite[Lemma~3.3]{Markering2026}. Each factor has a one-third logarithmic correction, giving the two-thirds correction in the capacity scale. The limiting normalized mean of their product under the stationary law is the capacity coefficient:
\begin{equation}
\kappa:=\mathsf E_{\mathsf Q}[X_\infty X_\infty^+]. \label{min:kappa}
\end{equation}
Here $X_\infty$ and $X_\infty^+$ are the normalized limits of the probabilities that a walk from $0$ avoids $\eta$ and $\eta[1,\infty)$, respectively. The constant $\alpha$ determines the conversion from walk steps to LERW edge count, while $\kappa$ is the capacity coefficient per unit of that length.

\begin{corollary}[Capacity limit and normalization]\label{min:capacity-prop}
The coefficient $\kappa$ is finite and positive, and
\begin{equation}\label{min:escape-constant-check}
\alpha\kappa =\mathsf E_{\mathsf P}[X_\infty^2X_\infty^+] =\frac{\pi^2}{12}.
\end{equation}
For every $t>0$,
\begin{equation}
\frac{L_n^2}{n}\operatorname{cap}(\eta[0,\lfloor nt\rfloor]) \longrightarrow\kappa t \quad\hbox{in probability under }\mathsf P. \label{min:capacitylaw}
\end{equation}
For the stationary two-sided path, the convergence holds in $L^1(\mathsf Q)$.  In particular,
\begin{equation}
M_N^{\mathrm{cap}}:=\mathsf E_{\mathsf Q}\operatorname{cap}(\widehat\eta[0,N]) \sim\kappa\frac{N}{(\log N)^{2/3}}. \label{min:expectedcap}
\end{equation}
\end{corollary}
\begin{proof}
\emph{Normalization of the capacity constant.} We first verify that $X_m^+\to X_\infty^+$ in every finite $L^p(\mathsf P)$. Decomposition at the last visit to $0$ gives
\begin{equation}\label{min:raw-green-factor}
\mathrm{Es}_m^+=\sum_{k=0}^m v_k\mathrm{Es}_{m-k},\qquad \mathrm{Es}_m^+\le G(0,0)\mathrm{Es}_{\lfloor m/2\rfloor}+C/m.
\end{equation}
In the first identity, the terms with $k\le m/2$, multiplied by $(\log m)^{1/3}$, converge almost surely to $G_\eta X_\infty$: use $X_m\to X_\infty$, summability of $(v_k)$ and dominated convergence. The remaining terms are at most $\sum_{k>m/2}p_k(0,0)\le C/m$. The second bound and the moment limits in \eqref{prelim:raw-escape-normalization} give uniform integrability of every fixed power, proving $X_m^+\to X_\infty^+$ in every finite $L^p(\mathsf P)$.

The density \eqref{min:density}, the moments of $X_\infty$ and $G_\eta\le G(0,0)$ give finiteness of $\kappa$, and positivity follows from $X_\infty^+\ge X_\infty>0$ $\mathsf Q$-almost surely. The density \eqref{min:density} and the moment limits with the testing horizon measured in random-walk steps \eqref{prelim:raw-escape-normalization} and \eqref{min:raw-green-factor} give
\[
\alpha\kappa=\mathsf E_{\mathsf P}[X_\infty^2X_\infty^+]=\lim_{m\to\infty}\mathsf E_{\mathsf P}[X_m^2X_m^+].
\]
Conditional on $\eta$, the last product is $(\log m)$ times the probability that two independent testing walks avoid $\eta$ through step $m$ and a third avoids $\eta\setminus\{0\}$. The square-horizon normalization in \cite[Theorem~1.2]{LawlerSunWu2019} has limit $\pi^2/24$. Each of the three escape factors changes by $2^{1/3}$ under our normalization by random-walk time, as explained after \eqref{prelim:raw-escape-normalization}. Their product therefore gains a factor of $2$, proving \eqref{min:escape-constant-check}.

\emph{Capacity convergence.} By \cite[Lemmas~3.4 and~3.6]{Markering2026}, the escape probabilities from $\widehat\eta[0,N]$ and $\widehat\eta[1,N]$ for a walk started at $0$, multiplied by $(\log N)^{1/3}$, converge $\mathsf Q$-almost surely to $X_\infty$ and $X_\infty^+$, respectively. The capacity law \cite[Theorem~1.8]{Markering2026} therefore gives
\[
C_N(\eta):=\frac{(\log N)^{2/3}}{N}\operatorname{cap}(\eta[0,N])\longrightarrow\kappa\qquad\mathsf P\text{-almost surely}.
\]
Since $\mathsf Q\ll\mathsf P$, the same limit holds $\mathsf Q$-almost surely. Moreover, \cite[equation~(3.41)]{Markering2026} gives $\limsup_N\mathsf E_{\mathsf Q}C_N\le\kappa$, while Fatou's lemma gives the reverse lower bound. Hence $\mathsf E_{\mathsf Q}C_N\to\kappa$, proving \eqref{min:expectedcap}, and
\[
\mathsf E_{\mathsf Q}|C_N-\kappa|=\mathsf E_{\mathsf Q}C_N+\kappa-2\mathsf E_{\mathsf Q}(C_N\wedge\kappa)\longrightarrow0.
\]
Finally, substitute $N=\lfloor nt\rfloor$ and use $N/n\to t$ and $\log(n+2)/\log N\to1$. This proves \eqref{min:capacitylaw} and its stationary $L^1$ version.
\end{proof}

\subsection{From whole-path escape to local segments}

A short window allows more escape than the whole finite path. The next proposition shows that the normalized sum of these excesses vanishes.

\begin{proposition}[Local-to-global replacement]\label{min:scalar}
For each $n\ge1$, let $c_n$ be as in \eqref{not:cutoffs}. For a one-sided simple path $\gamma$ and an integer $m\ge0$, define
\[
\Delta_{m,c_n}(\gamma)=\sum_{j=c_n}^{m-c_n} \left[e_{\gamma[j-c_n,j+c_n]}(\gamma_j) -e_{\gamma[0,m]}(\gamma_j)\right],
\]
where the sum is zero when $m<2c_n$. Then, for every $0<T_0<\infty$,
\begin{equation}
\sup_{0\le t\le T_0}\frac{L_n^2}{n} \Delta_{\lfloor nt\rfloor,c_n}(\eta)\longrightarrow0 \quad\hbox{in }\mathsf P\hbox{-probability}. \label{min:scalar-eq}
\end{equation}
\end{proposition}
\begin{proof}
We first prove the estimate at the largest length under $\mathsf Q$, transfer it to $\mathsf P$, and then use monotonicity in the length to obtain the supremum. Put $m_n=\lfloor nT_0\rfloor$ and use $c_n$ from \eqref{not:cutoffs}. Under the stationary law, put $a(k,\ell)=\mathsf E_{\mathsf Q}e_{\widehat\eta[-k,\ell]}(0)$. The equilibrium-measure identity and stationarity give $M_N^{\mathrm{cap}}=\sum_{j=0}^Na(j,N-j)\ge(N+1)a(N,N)$, by monotonicity of avoidance. Every summand of $\Delta_{m_n,c_n}$ is nonnegative, and the omitted boundary equilibrium masses total at most $2c_n$, so
\[
0\le\mathsf E_{\mathsf Q}\Delta_{m_n,c_n} \le(m_n-2c_n+1)a(c_n,c_n)-M_{m_n}^{\mathrm{cap}}+2c_n.
\]
The bound $a(c_n,c_n)\le M_{c_n}^{\mathrm{cap}}/(c_n+1)$ and \eqref{min:expectedcap}, together with $\log c_n/\log n\to1$ and $c_nL_n^2/n\to0$, make the normalized right-hand side tend to zero. Thus $L_n^2\Delta_{m_n,c_n}/n\to0$ in $L^1(\mathsf Q)$.

To transfer the result to $\mathsf P$, fix $k$. The local-window sums for $\gamma$ and $\vartheta^k\gamma$ have at most $2k$ unmatched summands, each bounded by one. Adding the omitted at most $2c_n$ boundary equilibrium terms to each global sum gives the capacities of $\gamma[0,m_n]$ and $\gamma[k,k+m_n]$, which differ by at most $2k$. Consequently
\[
|\Delta_{m_n,c_n}(\vartheta^k\gamma)-\Delta_{m_n,c_n}(\gamma)|\le4k+4c_n.
\]
After normalization the shift error vanishes, and Lemma~\ref{min:support} yields convergence in $\mathsf P$-probability. Finally, for each fixed $n$, the function $t\mapsto\Delta_{\lfloor nt\rfloor,c_n}(\gamma)$ is nondecreasing: extending $\gamma[0,\lfloor nt\rfloor]$ decreases the old whole-path escape probabilities and adds nonnegative summands. The supremum in \eqref{min:scalar-eq} is therefore its value at $T_0$.
\end{proof}

We now combine the capacity law and the replacement estimate to prove Proposition~\ref{bridge:row-convergence}.

\begin{proof}[Proof of Proposition~\ref{bridge:row-convergence}]
We first prove convergence of cumulative attachment masses, then pass from infinite-time local escape to early avoidance. Write $e_j=e_{\mathcal C_j^\eta}(\eta_j)$ and $m_{n,t}=\lfloor nt\rfloor$. For $m\ge2c_n$, the interior local sum differs from the interior equilibrium sum by $\Delta_{m,c_n}$, and each of the two omitted boundary sums is at most $2c_n$. For $m<2c_n$, both complete sums have at most $2c_n$ terms. Thus, in either case,
\[
\left|\frac{L_n^2}{n}\sum_{j=0}^m e_j-\frac{L_n^2}{n}\operatorname{cap}(\eta[0,m])\right|\le\frac{L_n^2}{n}\Delta_{m,c_n}+\frac{4c_nL_n^2}{n}.
\]
Proposition~\ref{min:scalar}, \eqref{min:capacitylaw} and $c_nL_n^2/n\to0$ give $(L_n^2/n)\sum_{j\le m_{n,t}}e_j\to\kappa t$ in probability for each $t>0$; at $t=0$ the sum is at most $L_n^2/n$. These cumulative sums are nondecreasing. A finite grid and continuity of $t\mapsto\kappa t$ therefore make the convergence uniform on $[0,T]$.

The deterministic bound \eqref{bridge:earlyrowerror} controls the total change from $e_j$ to $p_j$. Hence their cumulative masses also converge uniformly, which proves the weak convergence in \eqref{bridge:rowmeasure}. Finally, for $1\le p<3/2$, choose $p<p'<3/2$. The $L^{p'}$ bound \eqref{bridge:eqrowLp} makes the $p$th powers of the total masses uniformly integrable, so convergence in probability implies convergence in $L^p$.
\end{proof}

\section{The ordinary initial-spine insertion measure}
\label{sec:bridge}
Proposition~\ref{bridge:row-convergence} identifies the limiting attachment measure. We now record the LERW lifetime and the entire spatial path of the inserted branch. The main proposition below gives the insertion measure when no earlier branches have been exposed, together with the conditional-volume estimate needed for the vertex-retention rule in Section~\ref{sec:closure}. These conclusions and the ordinary spine invariance principle initiate the finite-history argument there. The random-walk time bounds proved here also apply to later insertions.

Fix $0<T<\infty$ and set $q_n=\lfloor nT\rfloor$. We use the local windows $\mathcal C_j^\eta$, escape probabilities $e_j,p_j$, and normalized sums $Z_{n,T},E_{n,T}$ defined in Section~\ref{prelim:lerw}. Throughout this section, testing-walk probabilities and expectations are conditional on the spine $\eta$. Let $W$ start at $\eta_j$ and put
\[
J_k=|\operatorname{LE}(W[0,k])|,\qquad \tau_\eta^+=\inf\{k\ge1:W_k\in\eta[0,\infty)\}.
\]
Write $\Gamma^{0}_{n,j,k}$ for its loop-erased path translated by $-\eta_j$, with space divided by $R_n$ and the LERW length parameter by $n$, linearly interpolated and held constant after its lifetime $J_k/n$. The lifetime is included in the recorded stopped path, also called its mark. We use the stopped-path topology of Section~\ref{prelim:lerw}. In the proposition, $B$ is Brownian motion with covariance $\Sigma=I_4/4$, $B^{[t]}$ has lifetime $t$ and is held constant after $t$, and $d_{\rm BL}$ is a bounded-Lipschitz metric for weak convergence of finite measures on attachment coordinates and stopped paths.

\begin{proposition}[The marked initial-spine insertion measure]\label{bridge:initial-kernel}
Define the finite random measure on attachment coordinates and relative stopped paths by
\[
I_{n,T}^{\eta} =b_n\sum_{j=0}^{q_n}\sum_{k\ge1} \Ee_{\eta_j}\!\left[ \one_{\{k<\tau_\eta^+,\ j+J_k\le q_n\}} \delta_{(j/n,\Gamma^{0}_{n,j,k})}\right].
\]
Then, for every $1\le p<3/2$,
\begin{equation}
d_{\rm BL}(I_{n,T}^{\eta},I_T)\longrightarrow0 \quad\hbox{in }L^p, \qquad I_T=\kappa\int_0^T\dd s\int_0^{T-s}\dd t\, \delta_s\otimes\Law(B^{[t]}). \label{bridge:fullkernel}
\end{equation}
Jointly with this convergence, the linearly interpolated rescaled spine defined by $\xi_n(j/n)=R_n^{-1}\eta_j$, $j\ge0$, converges locally uniformly in distribution to a Brownian spine with covariance $\Sigma=I_4/4$. In particular the translation-covariant conditional first moment supplied by Wilson's formula satisfies
\begin{equation}
Y_{n,T}(\eta):=\Ee[b_n|B_{\mathcal U}(0,q_n)|\mid\eta] =b_n(q_n+1)+\|I_{n,T}^{\eta}\| \longrightarrow\kappa T^2/2 \quad\hbox{in }L^p,\quad 1\le p<3/2. \label{bridge:firstmoment}
\end{equation}
\end{proposition}

The triangular region $s+t\le T$ in the limit imposes the intrinsic-ball restriction: the distance $s$ from the root to the attachment plus the new branch length $t$ must be at most $T$. Note that $I_T$ is the continuum insertion measure \eqref{min:continuum-row} for the spine alone, written in relative coordinates.

On compact intervals of walk time, we replace avoidance of the spine by early local avoidance (Lemma~\ref{bridge:microTV}) and identify the shape of the branch (Lemma~\ref{bridge:weightedFCLT}); together with Proposition~\ref{bridge:row-convergence} this gives the limit on such intervals (Lemma~\ref{bridge:compact-marked-kernel}). We then show that very short and very long walks contribute little (Lemma~\ref{bridge:raw-cutoffs}). Finally, we impose the intrinsic-ball restriction and remove the walk-time cutoffs to prove Proposition~\ref{bridge:initial-kernel}.

\subsection{A local comparison on compact random-walk time intervals}
\label{bridge:microscopic-section}

We first compare avoidance of the whole spine with avoidance of a local window during the first $r_n$ walk steps. We then show that the initial path segment $W[0,r_n]$ has negligible effect on the rescaled stopped loop-erasure. Together with \eqref{bridge:rowmeasure}, these comparisons determine the insertion measure on compact intervals of random-walk time.

Fix $0<\varepsilon<M<\infty$ and set
\[
\mathcal I_n=\{\lceil\varepsilon N_n\rceil,\ldots,\lfloor MN_n\rfloor\},\qquad E_j=\{W[1,r_n]\cap \mathcal C_j^\eta=\varnothing\},\qquad G_{j,k}=\{k<\tau_\eta^+\}.
\]
Let $\Omega_{\rm walk}$ be the countable space of finite nearest-neighbour paths in $\Z^4$, with arbitrary starting point and its discrete $\sigma$-field. For every bounded measurable $F$ on $\mathsf X=[0,T]\times[\varepsilon,M]\times\Omega_{\rm walk}$, define
\[
\begin{aligned}
\mathcal O_n^{\varepsilon,M}(F) &=b_n\sum_{j=0}^{q_n}\sum_{k\in \mathcal I_n} \Ee_{\eta_j}\!\left[ \one_{G_{j,k}}F(j/n,k/N_n,W[0,k])\right],\\
\mathcal O_n^{\rm loc,\varepsilon,M}(F) &=b_n\sum_{j=0}^{q_n}\sum_{k\in \mathcal I_n} \Ee_{\eta_j}\!\left[ \one_{E_j}F(j/n,k/N_n,W[0,k])\right].
\end{aligned}
\]
The two measures record the attachment, rescaled number of walk steps and complete testing path. They differ only in the avoidance requirement: $G_{j,k}$ requires avoidance of the whole spine through step $k$, whereas $E_j$ requires avoidance of the local window through step $r_n$. The intrinsic-ball restriction $j+J_k\le q_n$ is imposed in Section~\ref{bridge:initial-kernel-section}.

\begin{lemma}\label{bridge:microTV}
For fixed $T,\varepsilon,M$,
\begin{equation}
\|\mathcal O_n^{\rm loc,\varepsilon,M} -\mathcal O_n^{\varepsilon,M}\|_{\rm TV}\longrightarrow0 \quad\hbox{in probability and in }L^1. \label{bridge:microTVeq}
\end{equation}
The same statement holds after any measurable marking of the random-walk path, including its complete finite loop-erasure.
\end{lemma}
\begin{proof}
We compare the nonnegative difference of the two occupation measures, first for a finite initial spine segment and then for the remaining ray. Let $\Delta_n$ be the distance in \eqref{bridge:microTVeq}. For sufficiently large $n$, every $k\in\mathcal I_n$ satisfies $k\ge r_n$. Since $G_{j,k}\subset E_j$, their difference is nonnegative and, with $\|\mu\|_{\rm TV}=\sup_{|F|\le1}|\mu(F)|$,
\[
\Delta_n =b_n\sum_{j=0}^{q_n}\sum_{k\in \mathcal I_n} \Pp_{\eta_j}(E_j,\tau_\eta^+\le k) =b_n\sum_{j=0}^{q_n}\sum_{k\in \mathcal I_n} [p_j-\Pp_{\eta_j}(G_{j,k})].
\]
\emph{The finite initial segment.} Truncate the spine at a fixed LERW length $An$, with $A>T+1$. Put $m_{n,A}=\lfloor An\rfloor$, $\mathcal S_{n,A}=\eta[0,m_{n,A}]$ and $\mathcal R_{n,A}=\eta[m_{n,A}+1,\infty)$, viewed as vertex sets. Write $g_{A,j}(k)=\Pp_{\eta_j}(\tau_{\mathcal S_{n,A}}^+>k)$. For large $n$, $q_n\le m_{n,A}-c_n$, $k\ge r_n$ on $\mathcal I_n$, and
\[
\Pp_{\eta_j}(\tau_\eta^+>k)\le g_{A,j}(k)\le p_j.
\]
Thus $\Delta_n=\Delta_{n,A}^{\rm fin}+\Delta_{n,A}^{\rm tail}$, where
\[
\begin{aligned}
\Delta_{n,A}^{\rm fin} &=b_n\sum_{j=0}^{q_n}\sum_{k\in \mathcal I_n}[p_j-g_{A,j}(k)],\\
\Delta_{n,A}^{\rm tail} &=b_n\sum_{j=0}^{q_n}\sum_{k\in \mathcal I_n} [g_{A,j}(k)-\Pp_{\eta_j}(\tau_\eta^+>k)].
\end{aligned}
\]
The first term compares local-window avoidance with finite-spine avoidance; the second compares finite-spine avoidance with whole-spine avoidance. Both are nonnegative. For $c_n\le j\le q_n$, the first term is controlled by
\[
p_j-g_{A,j}(k) \le(p_j-e_j)+(e_j-e_{\mathcal S_{n,A}}(\eta_j)).
\]
The second differences form a sub-sum of $\Delta_{m_{n,A},c_n}(\eta)$. The at most $c_n$ indices $j<c_n$ contribute at most $b_nc_n|\mathcal I_n|$. Since $|\mathcal I_n|/N_n\le M$, $b_nN_n=L_n^2/n$ and $c_nL_n^2/n=O(L_n^{-38})$, \eqref{bridge:earlyrowerror} and Proposition~\ref{min:scalar} give
\[
\Delta_{n,A}^{\rm fin} \le M\left[ E_{n,T}-Z_{n,T}+\frac{L_n^2}{n}\Delta_{m_{n,A},c_n}(\eta) +\frac{c_nL_n^2}{n}\right]\longrightarrow0 \quad\hbox{in probability}.
\]
This proves the finite-spine comparison for each fixed $A$.

\emph{The remaining ray.} To restore the tail, we first show that it is unlikely to return to a fixed spatial ball. For fixed $A,D>0$,
\begin{equation}
\limsup_{n\to\infty}\Pp\{ \eta[\lfloor An\rfloor,\infty)\cap B(0,DR_n)\ne\varnothing\} \le CD^2/A. \label{bridge:oldtail}
\end{equation}
Set $t_{n,A}=\lfloor AN_n/2\rfloor$ and $B=B(0,DR_n)\cap\Z^4$. The estimate for $\sigma_{m_{n,A}}$ in \eqref{bridge:ordinaryclock} gives $\Pp(\sigma_{m_{n,A}}<t_{n,A})=o(1)$, and otherwise the indicated LERW tail is contained in $S[t_{n,A},\infty)$. The standard last-exit, heat-kernel and ball-capacity bounds \cite[Proposition~4.6.4, third assertion; Theorem~2.1.3, equation~(2.8); Proposition~6.5.2]{LawlerLimic2010} give
\[
\Pp\{S[t_{n,A},\infty)\cap B\ne\varnothing\} \le C\operatorname{cap}(B)/t_{n,A}\le CD^2/A,
\]
which proves \eqref{bridge:oldtail}.

Choose radii $D>D_0+2$ and define
\[
U=\left\{\max_{j\le q_n}|\eta_j|\le D_0R_n\right\},\qquad V=\{\mathcal R_{n,A}\cap B(0,DR_n)=\varnothing\}.
\]
The inclusion $\eta[0,q_n]\subset S[0,\sigma_{q_n}]$, \eqref{bridge:ordinaryclock} applied to $\sigma_{q_n}$ with the walk-step threshold $\lceil2(T+1)N_n\rceil$, and \eqref{bridge:gaussian-maximal} give
\[
\limsup_n\Pp(U^c)\le Ce^{-D_0^2/[C(T+1)]},\qquad \limsup_n\Pp(V^c)\le CD^2/A.
\]
For a spine in $U\cap V$, let $H_j=\{\max_{\ell\le r_n}|W_\ell-\eta_j|\le R_n\}$. Full avoidance of $\mathcal S_{n,A}$ through step $k$ implies $E_j$. The increments after step $r_n$ are independent of $E_j\cap H_j$, whose probability is at most $p_j$. Applying \eqref{bridge:gaussian-maximal} to the first $r_n$ steps and to the subsequent increments gives
\[
\begin{aligned}
g_{A,j}(k)-\Pp_{\eta_j}(\tau_\eta^+>k) &=\Pp_{\eta_j}\{\tau_{\mathcal S_{n,A}}^+>k,\, W[1,k]\cap\mathcal R_{n,A}\ne\varnothing\}\\
&\le \Pp_{\eta_j}(H_j^c) +\Pp_{\eta_j}\!\left( E_j\cap H_j,\, \max_{\ell\le k-r_n}|W_{r_n+\ell}-W_{r_n}| \ge(D-D_0-1)R_n\right)\\
&\le Ce^{-R_n^2/(Cr_n)}+Cp_je^{-(D-D_0-1)^2/(CM)}.
\end{aligned}
\]
Summing over $j,k$, with $R_n^2/r_n\asymp L_n^{13}$, gives on $U\cap V$
\[
\Delta_{n,A}^{\rm tail} \le C_{T,M}L_n^2e^{-L_n^{13}/C} +C_M E_{n,T}e^{-(D-D_0-1)^2/(CM)}.
\]
\emph{Convergence in probability and in $L^1$.} The constants above do not depend on  $n,A,D,D_0$ or  the spine. Lemma~\ref{bridge:rows} makes $E_{n,T}$ tight. Given a tolerance, choose a bound for $E_{n,T}$ and then $D_0$ so that $U^c$ is unlikely; next choose $D$ so that the second term is small; then choose $A$ so that $V^c$ is unlikely; finally let $n\to\infty$. Together with the finite-spine comparison, this proves convergence in probability. Since $0\le\Delta_n\le ME_{n,T}$, the higher moment estimate \eqref{bridge:eqrowLp} makes this family uniformly integrable and gives $L^1$ convergence. Total variation contracts under measurable pushforward, proving the marked assertion.
\end{proof}

We now compare the complete loop-erasure, conditional on the initial path segment $W[0,r_n]$, with the loop-erasure of the independent walk after step $r_n$. Keeping the factor $p_j$ in the main comparison error makes the estimate summable over attachment points.

\begin{lemma}[Weighted comparison after the initial path segment]
\label{bridge:weightedFCLT}
Fix $0<T<\infty$ and $0<\varepsilon<M<\infty$, with $\mathcal I_n$ as above. Let $x\in\mathbb Z^4$, let $W$ start at $x$, and let $E\in\mathcal F_{r_n}^W$, with $p=\Pp_x(E)$. Let $\Gamma_{n,k}^{x,0}$ be the complete relative stopped loop-erasure of $W[0,k]$, translated by $-x$ and scaled by $R_n$ in space and $n$ in the LERW length parameter, including its lifetime $J_k/n$. There is a deterministic sequence $\epsilon_n^{\mathrm{path}}\to0$ and a constant $C$ such that, for every $s\in[0,T]$ and every bounded Lipschitz test $F(s,t,\beta)$ with supremum norm and Lipschitz constant at most one for the bounded product metric on attachment coordinates, random-walk times and stopped paths,
\begin{equation}\label{bridge:weighted-prefix}
\left|\Ee_x\!\left[\one_E F(s,k/N_n,\Gamma_{n,k}^{x,0})\right]-p\Ee\!\left[F(s,k/N_n,B^{[k/N_n]})\right]\right|\le p\epsilon_n^{\mathrm{path}}+Ce^{-L_n^{11}/C},\qquad k\in\mathcal I_n.
\end{equation}
The estimate is uniform in $x,s,E$, the tests and the indicated random-walk times. The same assertion holds conditionally on exposed information under which $W$ is a fresh simple random walk from $x$. In particular, conditionally on the spine, take $x=\eta_j$, $s=j/n$, and $E=E_j=\{W[1,r_n]\cap \mathcal C_j^\eta=\varnothing\}$; then $p=p_j$ and $\Gamma^0_{n,j,k}=\Gamma_{n,k}^{\eta_j,0}$.
\end{lemma}
\begin{proof}
Fix $k\in\mathcal I_n$. For sufficiently large $n$, $k\ge r_n+2u_n$. We first compare the complete loop-erasure with the loop-erasure of the independent walk after step $r_n$, and then integrate the comparison over the early event $E$ while retaining its probability $p$.

\emph{Matching the loop-erasures after a cut.} Condition on $\mathcal F_{r_n}^W$. The remaining increments form a fresh walk $S_l=W_{r_n+l}-W_{r_n}$, independent of $\mathcal F_{r_n}^W$. Define
\[
\begin{aligned}
\mathcal E_{n,1} &=\{(W_{r_n}+S[u_n,\infty))\cap W[0,r_n]=\varnothing\},\\
\mathcal E_{n,2} &=\{\text{$S$ has a global cut in $[u_n,2u_n]$}\}.
\end{aligned}
\]
The heat-kernel bound $\Pp(S_l=y)\le C(l+1)^{-2}$ and a union bound over the at most $r_n+1$ vertices of $W[0,r_n]$ give, almost surely,
\[
\Pp\bigl((W_{r_n}+S[u_n,\infty))\cap W[0,r_n]\ne\varnothing \mid\mathcal F_{r_n}^W\bigr)\le C(r_n+1)/u_n.
\]
Together with \eqref{bridge:cutbound}, this gives, almost surely,
\begin{equation}\label{bridge:join-failure}
\Pp\!\left(\mathcal E_{n,1}^{\mathrm c}\cup\mathcal E_{n,2}^{\mathrm c} \,\middle|\,\mathcal F_{r_n}^W\right) \le a_n^{\mathrm{join}} :=C(r_n+1)/u_n+C\log\log n/\log n=o(1).
\end{equation}
On $\mathcal E_{n,1}\cap\mathcal E_{n,2}$, choose a global cut $v\in[u_n,2u_n]$. Then $r_n+v$ is a cut of the concatenated path. Hence $\operatorname{LE}(W[0,k])$ and $W_{r_n}+\operatorname{LE}(S[0,k-r_n])$ have exactly the same ordered suffix after this cut, while their preceding portions have at most $r_n+4u_n$ edges in total. Therefore, their lifetimes and the rescaled LERW length parameters of matched vertices differ by at most $(r_n+4u_n)/n=o(1)$.

\emph{Control in the stopped-path metric.} It remains to control the unmatched initial portions and the small time shift of the common suffix. Take $\theta_n=L_n^{-1}$. By \eqref{bridge:gaussian-maximal},
\[
\Pp_x\!\left\{\max_{l\le r_n}|W_l-x|>\theta_nR_n\right\}\le Ce^{-L_n^{11}/C},\qquad \Pp\!\left\{\max_{l\le2u_n}|S_l|>\theta_nR_n\right\}\le Ce^{-L_n^5/C}.
\]
Outside these two exceptional events, both the diameters of the unmatched portions and the magnitude of the translation $W_{r_n}-x$ are $O(\theta_nR_n)$, and matched vertices have a shift of at most $(r_n+4u_n)/n=o(1)$ in the LERW length parameter. Corollary~\ref{bridge:finiteFCLT} gives weak convergence of the fresh stopped-path laws, uniformly over durations $(k-r_n)/N_n\in[\varepsilon/2,M]$. The resulting tightness controls the modulus of continuity at this shift, and $r_n/N_n\to0$ permits replacing the fresh duration $(k-r_n)/N_n$ by $k/N_n$.

\emph{Integration over the early event.} Put $D=\{\max_{l\le r_n}|W_l-x|\le\theta_nR_n\}$. The preceding estimates are uniform over the initial paths in $D$: the comparison-failure probability is bounded by \eqref{bridge:join-failure}, and all remaining random errors involve the fresh walk. Together with Corollary~\ref{bridge:finiteFCLT}, they give a deterministic $\epsilon_n^{\mathrm{path}}\to0$ such that, on $D$,
\[
\left|\Ee_x\!\left[F(s,k/N_n,\Gamma_{n,k}^{x,0})\mid\mathcal F_{r_n}^W\right]-\Ee F(s,k/N_n,B^{[k/N_n]})\right|\le\epsilon_n^{\mathrm{path}},
\]
uniformly over $k$ and the stated tests. Since $E\in\mathcal F_{r_n}^W$, multiplying by $\mathbf1_E$ and integrating over $D$ costs at most $p\epsilon_n^{\mathrm{path}}$. On $D^{\mathrm c}$ the tests are bounded by one, so the remaining error is at most $2\Pp_x(D^{\mathrm c})$. The first Gaussian estimate therefore proves \eqref{bridge:weighted-prefix}. Translation invariance makes the estimate uniform in $x$, and the same argument applies under the stated conditioning.
\end{proof}

Combining the two comparisons with the attachment measure yields the insertion measure on compact random-walk time intervals.

\begin{lemma}[Marked occupation on compact random-walk time intervals]
\label{bridge:compact-marked-kernel}
Fix $0<T<\infty$ and $0<\varepsilon<M<\infty$. Let $\mathcal O_{n,\mathrm{LE}}^{\rm loc,\varepsilon,M}$ and $\mathcal O_{n,\mathrm{LE}}^{\varepsilon,M}$ be the pushforwards of the corresponding random-walk path measures under $W[0,k]\mapsto\Gamma^{0}_{n,j,k}$, retaining attachment and random-walk time. As finite measures on attachment coordinates, random-walk times and relative stopped paths, both converge in probability and in $L^1$ for $d_{\rm BL}$ to
\[
\mathcal O^{\varepsilon,M} =\kappa\int_0^T\!\dd s\int_\varepsilon^M\!\dd t\, \delta_{(s,t)}\otimes\Law(B^{[t]}), \qquad \Ee B_tB_t^{\mathsf T}=tI_4/4.
\]
In particular, for every bounded continuous stopped-path test $F$,
\begin{equation}\label{bridge:compactkernel}
b_n\sum_{j=0}^{q_n}\sum_{k\in\mathcal I_n}\Ee_{\eta_j}\!\left[\one_{G_{j,k}}F(j/n,k/N_n,\Gamma^0_{n,j,k})\right]\xrightarrow{\mathbb P}\kappa\int_0^T\int_{\varepsilon}^M\Ee[F(s,t,B^{[t]})]\dd t\dd s.
\end{equation}
The same limit holds with $\one_{E_j}$ in place of $\one_{G_{j,k}}$. The convergences hold jointly with the ordinary rescaled spine on bounded LERW length-parameter intervals. Translating each relative path by the spatial label of its attachment gives the corresponding absolute-path limit.
\end{lemma}
\begin{proof}
We compare the local occupation measure with a measure having the same attachment weights and Brownian path marks, then replace those weights and the walk-time grid by their limits. Define
\[
\Lambda_n=b_n\sum_{j=0}^{q_n}\sum_{k\in\mathcal I_n}p_j\,\delta_{(j/n,k/N_n)}\otimes\Law(B^{[k/N_n]}).
\]
Multiplying \eqref{bridge:weighted-prefix} by $b_n$ and summing over $j,k$ gives
\[
d_{\rm BL}(\mathcal O_{n,\mathrm{LE}}^{\rm loc,\varepsilon,M},\Lambda_n)\le C_M\epsilon_n^{\mathrm{path}}E_{n,T}+C_{T,M}L_n^2e^{-L_n^{11}/C}\longrightarrow0\quad\hbox{in probability},
\]
because $E_{n,T}$ is tight by Lemma~\ref{bridge:rows}.

By $b_nN_n=L_n^2/n$, the attachment factor in $\Lambda_n$ is exactly $(L_n^2/n)\sum_jp_j\delta_{j/n}$, which converges to $\kappa\,ds$ by Proposition~\ref{bridge:row-convergence}; the time factor is $N_n^{-1}\sum_{k\in\mathcal I_n}\delta_{k/N_n}$. For unit bounded-Lipschitz $F$, the functions $(s,t)\mapsto\Ee F(s,t,B^{[t]})$ have a common modulus of continuity on $[0,T]\times[\varepsilon,M]$. A fixed finite rectangular partition reduces the comparison to attachment-interval masses and time-grid counts. Since the limit has no atoms, the former converge in probability by \eqref{bridge:rowmeasure}, and the latter are Riemann sums. On sets where $E_{n,T}$ is bounded, the common modulus controls both approximation errors uniformly in $F$. Tightness of $E_{n,T}$ and refinement of the partition therefore give $d_{\rm BL}(\Lambda_n,\mathcal O^{\varepsilon,M})\to0$ in probability.

Lemma~\ref{bridge:microTV}, followed by the common loop-erasure marking, transfers the result to the full-avoidance measure. The total masses of both occupation measures are at most $ME_{n,T}$, so \eqref{bridge:eqrowLp} makes their distances from the finite limiting measure uniformly integrable. Their convergence thus also holds in $L^1$.

The limiting relative-path measure is deterministic. Slutsky's theorem and the ordinary spine invariance principle therefore give the asserted joint convergence. Finally, for a continuous spine $\xi$, the map $(s,t,\beta)\mapsto(s,t,\xi(s)+\beta)$ is continuous, also when the spine varies locally uniformly. Applying this common translation to the two measures proves the absolute-path assertion.
\end{proof}

\subsection{Removing the random-walk time cutoffs}
\label{bridge:cutoff-section}

A walk producing a branch of bounded LERW length may spend many walk steps in erased loops. We control the total occupation outside compact random-walk time intervals using $H_{A,a}(J,N)$ from Corollary~\ref{closure:extreme}. The bounds are uniform in the killing set, so we state them for arbitrary attachment sites and reuse them after a finite Wilson history has been exposed.

\begin{lemma}[Bounds for time cutoffs and omitted attachment indices]
\label{bridge:raw-cutoffs}
Fix $0<B<\infty$. Let $\mathcal J_n$ contain at most $Bn$ indices. For each $j\in\mathcal J_n$, fix a starting site $a_j$ and vertex sets $a_j\in \mathcal C_j\subseteq A_j\subseteq\Z^4$, and define
\[
p_j=\Pp_{a_j}(W[1,r_n]\cap \mathcal C_j=\varnothing),\qquad E_n=\frac{L_n^2}{n}\sum_{j\in\mathcal J_n}p_j.
\]
For an integer $1\le J\le Bn$ and a set of walk-step indices $\mathcal R\subseteq\mathbb N$, write
\[
\mathcal V_n(\mathcal R) =b_n\sum_{j\in\mathcal J_n}\sum_{k\in\mathcal R} \Pp_{a_j}\{k<\tau_{A_j}^+,\ 1\le J_k\le J\}, \qquad J_k=|\operatorname{LE}(W[0,k])|.
\]
There are constants $M_B,C_B<\infty$, independent of all sites and sets, such that, for all sufficiently large $n$ and all $\varepsilon>0$,
\begin{equation}
\mathcal V_n(\{k<\varepsilon N_n\}) \le\varepsilon E_n+C_BL_n^{-11}, \label{bridge:smallraw}
\end{equation}
and
\begin{equation}
\mathcal V_n(\{k>M_BN_n\}) \le C_B\frac{\log\log n}{L_n}E_n+C_B/L_n. \label{bridge:largeraw}
\end{equation}
Consequently $\mathcal V_n(\mathbb N)\le C_BE_n+C_B/L_n$. Deleting any $m_0$ attachment indices changes the entire normalized occupation by at most $C_Bm_0L_n^2/n$, even without local sets $\mathcal C_j$. All bounds hold conditional on any information determining the sites and sets, with independent testing walks. In particular, a uniform $L^p$ bound for $E_n$ implies the same bound for $\mathcal V_n(\mathbb N)$.
\end{lemma}
\begin{proof}
We prove the short-time bound first. For the long-time bound, we use a weighted pointwise estimate only up to order $n\log n$ and apply \eqref{bridge:universalH} to all later times.

\emph{Small walk times.} For $k<r_n$, the normalized occupation is at most $b_n|\mathcal J_n|r_n\le C_BL_n^{-11}$. For $r_n\le k<\varepsilon N_n$, avoidance of $A_j$ implies $W[1,r_n]\cap \mathcal C_j=\varnothing$. Summation and $b_nN_n=L_n^2/n$ prove \eqref{bridge:smallraw}.

\emph{Long walk times up to order $n\log n$.} Use the fresh walk $S_l=W_{r_n+l}-W_{r_n}$ and the events $\mathcal E_{n,1},\mathcal E_{n,2}$ defined in the proof of Lemma~\ref{bridge:weightedFCLT}. They respectively require avoidance of the initial path by the fresh tail and a global cut in $[u_n,2u_n]$. Given $\mathcal F_{r_n}^W$, \eqref{bridge:join-failure} bounds the conditional probability of $(\mathcal E_{n,1}\cap\mathcal E_{n,2})^{\mathrm c}$ by $C_B\log\log n/L_n^3$, almost surely. On $\mathcal E_{n,1}\cap\mathcal E_{n,2}$, $J_k\le J$ implies $|\operatorname{LE}(S[0,k-r_n])|\le J+r_n+4u_n$. Put $\ell_n=J+r_n+4u_n+1$. The part of the infinite loop-erasure retained by walk step $k-r_n$ is an initial segment of the finite loop-erasure, so
\[
\{J_k\le J\}\subseteq\{\sigma_{\ell_n}^{S}>k-r_n\} \quad\hbox{on }\mathcal E_{n,1}\cap\mathcal E_{n,2}.
\]
Since $4u_n\le\ell_n\le Bn+o(n)$, one has $\log\ell_n/\log n\to1$ uniformly in $J$. Choose a fixed $M_B$ so that $k-r_n>C_{\mathrm{time}}\ell_nL_{\ell_n}$ whenever $k>M_BN_n$. The estimate for $\sigma_{\ell_n}^S$ in \eqref{bridge:eqclock} bounds this probability by $C_B\log\log n/L_n^3$, uniformly in $J,k$. This bound and the comparison-failure bound hold conditionally on $\mathcal F_{r_n}^W$. Avoidance of $A_j$ through step $k$ requires the $\mathcal F_{r_n}^W$-measurable event $W[1,r_n]\cap \mathcal C_j=\varnothing$, of probability $p_j$. Multiplying the conditional bound by its indicator and taking expectation gives
\[
\Pp_{a_j}\{k<\tau_{A_j}^+,\ J_k\le J\} \le C_Bp_j\frac{\log\log n}{L_n^3}, \qquad k>M_BN_n.
\]
This pointwise bound is summed only over $M_BN_n<k\le A_Bn\log n$, with fixed $A_B$ to be chosen below. Since $(n\log n)/N_n=O(L_n^2)$, its contribution is at most
\[
C_B\frac{\log\log n}{L_n}E_n.
\]
\emph{Walk times beyond order $n\log n$.} Choose $A_B$ large enough for \eqref{bridge:universalH} to apply at $N=A_Bn\log n$ for every $J\le Bn$. That occupation bound controls each attachment's contribution from all remaining walk-step indices by $C_Bn/L_n^2$. Multiplication by $b_n|\mathcal J_n|$ gives at most $C_B/L_n$, proving \eqref{bridge:largeraw}.

\emph{Total occupation and omitted attachment indices.} The contribution from $r_n\le k\le M_BN_n$ is at most $M_BE_n$; together with the small- and large-time bounds, this proves the total bound. For any $m_0$ attachment indices, the first bound in \eqref{bridge:universalH}, with $J\le Bn$, gives $b_nm_0CJ L_J\le C_Bm_0L_n^2/n$, proving the deletion bound. The estimates hold for each fixed configuration; hence they also hold conditionally, and the total bound gives the moment assertion.
\end{proof}

\subsection{Proof of the initial-spine insertion limit}\label{bridge:initial-kernel-section}

\begin{proof}[Proof of Proposition~\ref{bridge:initial-kernel}]
\emph{Step 1: uniform moment bounds.} Apply Lemma~\ref{bridge:raw-cutoffs} with $B=T+1$, attachment sites $\eta_j$ for $0\le j\le q_n$, local sets $\mathcal C_j^\eta$, killing set $\eta[0,\infty)$ and length bound $J=q_n$. Then $E_n=E_{n,T}$. Denote the resulting occupation by $V_{n,T}^{\{k\in\mathcal R\}}=\mathcal V_n(\mathcal R)$ and $V_{n,T}=\mathcal V_n(\mathbb N)$. Its total bound and \eqref{bridge:eqrowLp} give $\sup_n\Ee V_{n,T}^p<\infty$ for every $1\le p<3/2$.

\emph{Step 2: imposing the intrinsic-ball restriction.} Fix $M>\max\{T,M_{T+1}\}$ and $0<\varepsilon<M$. On the observation space of Lemma~\ref{bridge:compact-marked-kernel}, write $v=\operatorname{lifetime}(\beta)$ and $\mathsf p(s,t,\beta)=(s,\beta)$. The part of $I_{n,T}^\eta$ with $\varepsilon N_n\le k\le MN_n$ is exactly
\[
\mathsf p_*\!\left(\mathbf1_{\{s+v\le q_n/n\}}\mathcal O_{n,\mathrm{LE}}^{\varepsilon,M}\right).
\]
The lifetime is a continuous coordinate. Moreover, under $\mathcal O^{\varepsilon,M}$ it equals the recorded Brownian duration $t$, so
\[
\mathcal O^{\varepsilon,M}\{s+v=T\}=\kappa\int_0^T\int_\varepsilon^M\mathbf1_{\{s+t=T\}}\,dt\,ds=0.
\]
For each fixed boundary width, bracket the length indicator by continuous cutoffs of $s+v$ near $T$. Lemma~\ref{bridge:compact-marked-kernel} applies to these continuous tests; $q_n/n\to T$ and the zero boundary mass allow the width to decrease to zero. Thus the displayed restricted measures converge in probability for $d_{\rm BL}$ to $\mathsf p_*(\mathbf1_{\{s+v\le T\}}\mathcal O^{\varepsilon,M})$, which is the part of $I_T$ with lifetime at least $\varepsilon$.

\emph{Step 3: removing the time cutoffs and proving $L^p$ convergence.} Apply \eqref{bridge:smallraw}--\eqref{bridge:largeraw} to the actual length-restricted discrete measure. They bound its omitted mass by
\[
\left(\varepsilon+C_T\frac{\log\log n}{L_n}\right)E_{n,T}+C_T(L_n^{-11}+L_n^{-1}).
\]
The omitted continuum mass below $\varepsilon$ is at most $\kappa T\varepsilon$, and there is no mass above $M>T$. Since $\sup_n\Ee E_{n,T}<\infty$, first $n\to\infty$ and then $\varepsilon\downarrow0$ give $I_{n,T}^\eta\to I_T$ in probability for $d_{\rm BL}$. For $1\le p<3/2$, choose $p<p'<3/2$ and use
\[
d_{\rm BL}(I_{n,T}^{\eta},I_T)\le\|I_{n,T}^{\eta}\|+\|I_T\|\le V_{n,T}+\kappa T^2/2.
\]
The $L^{p'}$ bound from Step 1 makes the $p$th powers uniformly integrable, proving \eqref{bridge:fullkernel}. Since $I_T$ is deterministic, Slutsky's theorem gives joint convergence with the ordinary rescaled spine.

Finally, the marked Wilson identity \eqref{min:occupation} partitions the root ball into its $q_n+1$ spine vertices and the off-spine insertions, giving the equality in \eqref{bridge:firstmoment}. The constant test transfers the $L^p$ convergence to total masses; $b_n(q_n+1)\to0$ and $\|I_T\|=\kappa T^2/2$ then prove the conditional first-moment assertion.
\end{proof}

\section{Convergence of the discrete Campbell measures}\label{sec:closure}
Section~\ref{sec:bridge} handles insertion along an ordinary infinite LERW. We extend that result to every sampled vertex by a mass-transport bound on its full future ray. Each added branch is a reversed prefix of such a ray, so the same insertion limit applies after restricting and changing coordinates. The remaining comparison concerns only a finite collection of other branches. Retaining vertices with bounded conditional ball volume makes all insertion masses uniformly bounded and allows us to iterate this argument.

Fix $T>0$ and $q_n=\lfloor nT\rfloor$. For each vertex $x$, let $\gamma^x$ be its full future ray in $\mathcal U$, translated to start at zero when used as an argument of a translation-covariant functional. For $t>0$, the marked Wilson identity specifies simultaneous versions of
\[
Y_{n,t}(\gamma^x) =\mathbb E\!\left[b_n|B_{\mathcal U}(x,\lfloor nt\rfloor)| \mid\gamma^x\right].
\]

\paragraph{Recording the exposed tree and sampled vertices.}\label{closure:history-coordinates}
After $k$ insertions, Wilson's algorithm has revealed the ray from $0$ to infinity and $k$ finite branches, each attached to an earlier one and ending at a sampled vertex. We record this information as a \emph{history}. Write $\gamma^0=\eta$ for the initial ray from $0$. For $1\le i\le k$, let $\gamma^i=(\gamma_0^i,\ldots,\gamma_{m_i}^i)$ be the finite simple path recorded at the $i$-th insertion, oriented from its attachment vertex in the previously exposed tree to the sampled vertex $\gamma_{m_i}^i=x_i$. If $x_i$ is already exposed, take $m_i=0$ and $\gamma^i=(x_i)$. Let $\pi_i\in\{0,\ldots,i-1\}$ be the smallest index such that $\gamma_0^i\in\gamma^{\pi_i}$, and let $j_i$ be the unique index with $\gamma_0^i=\gamma_{j_i}^{\pi_i}$. Thus $\pi_i$ specifies the parent path, $j_i$ specifies the attachment coordinate measured from that path's initial vertex, and $m_i$ is the number of new edges. The lattice history is $H=(\eta,(\pi_i,j_i,m_i,\gamma^i)_{i=1}^k)$, and its exposed tree is $\mathcal S(H)=\bigcup_{i=0}^k\gamma^i$. These records reconstruct the full future ray of each sample. Later attachments do not change the paths recorded at their creation.

At scaling index $n$, put $a_i=j_i/n$ and $\ell_i=m_i/n$. Define the spatial paths by $\xi(j/n)=R_n^{-1}\eta_j$ for $j\ge0$ and $\beta_i(j/n)=R_n^{-1}\gamma_j^i$ for $0\le j\le m_i$, using linear interpolation and holding $\beta_i$ constant after $\ell_i$. Write
\[
\widehat H_n=\operatorname{Scale}_n(H):=\left(\xi,(\pi_i,a_i,\ell_i,\beta_i)_{i=1}^k\right).
\]
A generic history in the rescaled coordinate space is denoted by $\mathsf h$. For either kind of rescaled record, put $\phi_0=\xi$, $\phi_i=\beta_i$ for $i\ge1$, and define $h_0=0$ and $h_i=h_{\pi_i}+a_i$. Construct its abstract tree by attaching the initial point of interval $[0,\ell_i]$ to coordinate $a_i$ of its parent path. Then $h_i$ is the root distance of the $i$-th attachment and $h_i+s$ is the root distance of coordinate $s$ on that branch. Consistency requires $\xi(0)=0$, $a_i\le\ell_{\pi_i}$ when $\pi_i>0$, $\beta_i(0)=\phi_{\pi_i}(a_i)$, and $h_i+\ell_i\le T$. Equal spatial labels do not make additional identifications in this abstract tree. Zero-length branches retain the sample order but add no tree points.

Let
\[
\mathscr P_T=\left\{(\ell,\beta)\in[0,T]\times C([0,\infty),\mathbb R^4):\beta(s)=\beta(s\wedge\ell)\text{ for every }s\ge0\right\}
\]
be the stopped-path space with the metric $d_{\rm stop}$ of Section~\ref{prelim:lerw}. The space of history coordinates is the Polish product
\[
\mathscr X_{k,T}=C([0,\infty),\mathbb R^4)\times\prod_{i=1}^k\left(\{0,\ldots,i-1\}\times[0,T]\times\mathscr P_T\right).
\]
Use locally uniform convergence for the initial ray, the discrete topology for parent indices, the usual topology for attachment coordinates, and the stopped-path topology for each $(\ell_i,\beta_i)$. Equip the product with the maximum of the corresponding bounded coordinate metrics. History measures are supported on the subset satisfying the preceding consistency conditions.  Appending an insertion retains its already chosen parent index and simply adds its attachment coordinate and stopped path.

For any $A>0$, let $\mathscr X_{k,T,A}$ be the product obtained from $\mathscr X_{k,T}$ by replacing the initial-ray space by $C([0,A],\mathbb R^4)$ with its uniform topology. For $\mathsf h=(\xi,(\pi_i,a_i,\ell_i,\beta_i)_{i=1}^k)$, define
\begin{equation}\label{closure:truncated-H}
\mathsf R_A(\mathsf h)=\left(\xi|_{[0,A]},(\pi_i,a_i,\ell_i,\beta_i)_{i=1}^k\right)\in\mathscr X_{k,T,A}.
\end{equation}
Thus $\mathsf R_A$ truncates only the initial ray and preserves every finite-branch coordinate. We use $A>T$ when these truncated coordinates must reconstruct a tree containing every recorded attachment.

For a lattice history $H$, let $\overline Q_{n,\mathcal S(H),T}$ be the complete conditional insertion measure in \eqref{min:wilson-row}, including atoms of mass $b_n$ at exposed vertices. The short length and time cutoffs $c_n,r_n,u_n$ are those of \eqref{not:cutoffs}. Throughout this section, estimate constants may depend on the fixed radius $T$, volume cutoff $K$ and number of samples $k$, but not on $n$ or the realized history. Any other dependence is stated explicitly.

\paragraph{Why restrict the vertex measure?}
The $k$-point Campbell measure has mass equal to the $k$th moment of normalized ball volume, whereas Proposition~\ref{bridge:initial-kernel} only gives conditional-volume control in $L^p$ for $p<3/2$. We obtain the bounds needed for iteration by retaining vertices whose conditional expected ball volume is at most a fixed $K$. This condition depends only on the vertex's future ray, so it is measurable as soon as that ray is exposed. Lemma~\ref{closure:good-deletion} will show that the volume removed in this way tends to zero. The radius $2T$ is used because its ball around any vertex of the root's radius-$T$ ball contains that entire root ball. Fix $K>2\kappa T^2$ and define
\[
\begin{gathered}
g_n(x)=\mathbf1_{\{Y_{n,2T}(\gamma^x)\le K\}},\qquad G_n=g_n(0),\qquad V_n=b_n|B_{\mathcal U}(0,q_n)|,\\
A_n=b_n\sum_{x\in B_{\mathcal U}(0,q_n)}g_n(x),\qquad D_n=V_n-A_n.
\end{gathered}
\]
A vertex $x$ is retained when $g_n(x)=1$: given its future ray, the conditional expected normalized volume of its radius-$2T$ ball is at most $K$. The retained measure on $B=B_{\mathcal U}(0,q_n)$ is
\begin{equation}
\widetilde\mu_{n,T} =G_n b_n\sum_{x\in B_{\mathcal U}(0,q_n)}g_n(x)\delta_x. \label{closure:eq3_2}
\end{equation}
For $k\ge0$, the discrete Campbell measure $\mathfrak C_{n,k,T}^{\,g}$ records the ordered retained vertices and their exposed paths. It is the discrete analogue of the measure $\mathfrak C_{k,T}$ of Corollary~\ref{min:continuum-campbell}, with the retention restriction imposed. It is the image of
\begin{equation}
\mathbb P(d\mathcal U)\,G_n\,b_n^k \sum_{x_1,\ldots,x_k\in B_{\mathcal U}(0,q_n)} \prod_{i=1}^kg_n(x_i)\, \delta_{(x_1,\ldots,x_k)} \label{closure:eq10_1}
\end{equation}
under the map recording the union of the full future rays and the ordered samples.

The main result of this section is the following convergence of the rescaled records of finitely many samples. Its proof is completed in Section~\ref{closure:induction}.

\begin{proposition}[Convergence of the retained discrete Campbell measures]
\label{closure:campbell-convergence}
Fix $T>0$ and $K>2\kappa T^2$. For every fixed integer $k\ge0$,
\begin{equation}
(\operatorname{Scale}_n)_*\mathfrak C_{n,k,T}^{\,g} \Longrightarrow\mathfrak C_{k,T}. \label{closure:campbell-limit}
\end{equation}
The convergence is weak as finite measures on $\mathscr X_{k,T}$. It records attachment indices and coordinates, lifetimes, entire stopped spatial paths, and the initial ray locally uniformly in the LERW length parameter.

Consequently, for every $k\ge1$, the expected rooted marked distance-matrix measures of the rescaled ball with vertex measure $\widetilde\mu_{n,T}$ converge to those of $B_T(\mathcal K_{\kappa,I_4/4})$.
\end{proposition}

Sections~\ref{a-truncation-with-vanishing-deleted-mass}--\ref{closure:screening} establish sampled-ray domination, identify the reference insertion kernels, and control the effect of the other exposed branches. Section~\ref{closure:induction} combines them into a one-step limit and proves Proposition~\ref{closure:campbell-convergence} by induction. The proof of Theorem~\ref{min:main} then applies the moment criterion of Lemma~\ref{prelim:marked-moments} and removes the vertex truncation.

\subsection{Retained vertices and bounded insertion measures}\label{a-truncation-with-vanishing-deleted-mass}

The retention condition in \eqref{closure:eq3_2} has two roles: it removes a vanishing expected volume and bounds the total mass of every subsequent insertion measure. The second property supplies the moment bounds needed for Campbell convergence. Both rely on the retention condition being measurable from the exposed future ray.

\begin{lemma}[Vanishing removed volume]\label{closure:good-deletion}
With the preceding definitions, $\Pp(G_n=0)\to0$ and
\begin{equation}
\delta_n:=\mathbb ED_n =\mathbb E[V_n(1-G_n)] =\mathbb E\!\left[ Y_{n,T}(\gamma^0) \mathbf1_{\{Y_{n,2T}(\gamma^0)>K\}}\right]\longrightarrow0, \label{closure:eq3_1}
\end{equation}
and
\begin{equation}
\mathbb E\!\left[V_n-\widetilde\mu_{n,T}(B)\right] =\mathbb E[V_n-G_nA_n]\le2\delta_n\longrightarrow0. \label{closure:eq3_3}
\end{equation}
\end{lemma}

\begin{proof}
Apply the mass-transport principle \cite[Definition~2.1 and Theorem~3.2]{AldousLyons2007}, as used in \cite[Lemma~6.5 and equation~(6.11)]{HutchcroftSousi2023}, to
\[
(x,y)\longmapsto b_n\mathbf1_{\{d_{\mathcal U}(x,y)\le q_n\}}(1-g_n(y)).
\]
It gives $\mathbb ED_n=\mathbb E[V_n(1-G_n)]$, and conditioning on $\gamma^0$ gives \eqref{closure:eq3_1}. Proposition~\ref{bridge:initial-kernel} gives $Y_{n,2T}\to2\kappa T^2<K$ in $L^1$ under the ordinary ray law. This implies both $\mathbb P(G_n=0)\to0$ and $\mathbb E[Y_{n,2T}\mathbf1_{\{Y_{n,2T}>K\}}]\to0$. The conditional-volume definition makes $Y_{n,t}$ nondecreasing in $t$, so $Y_{n,T}\le Y_{n,2T}$ and $\delta_n\to0$. Finally, $V_n-G_nA_n=D_n+(1-G_n)A_n\le D_n+(1-G_n)V_n$, and each term on the right has expectation $\delta_n$.
\end{proof}

We next use the retention condition to bound insertion masses and sampling moments.

\begin{lemma}[Bounded insertion measures and retained sampling moments]
\label{closure:bounded-rows}
\label{bounded-insertion-rows-and-all-moments-of-the-retained-measure}
For an ordered tuple $(x_0=0,x_1,\ldots,x_k)$, allowing repetitions, suppose that $x_i\in B_{\mathcal U}(0,q_n)$ and $g_n(x_i)=1$ for all $i$. Expose $\mathcal S_k=\bigcup_{i=0}^k\gamma^{x_i}$ by running Wilson's algorithm from these vertices first. Then the complete next insertion measure satisfies
\begin{equation}
\|\overline Q_{n,\mathcal S_k,T}\|=\mathbb E\!\left[b_n|B_{\mathcal U}(0,q_n)|\mid\mathcal S_k\right]\le\sum_{i=0}^kY_{n,2T}(\gamma^{x_i})\le(k+1)K.\label{closure:eq4_2}
\end{equation}
For every integer $k\ge0$ and every nonnegative measurable functional $F$ of the initial ray,
\begin{equation}
\mathbb E\!\left[G_nA_n^kF(\gamma^0)\right] \le k!K^k\,\mathbb E F(\gamma^0). \label{closure:eq4_3}
\end{equation}
In particular, $\|\mathfrak C_{n,k,T}^{\,g}\|\le k!K^k$. More generally, for every $0\le i\le k$ and every nonnegative measurable ray functional $F$,
\begin{equation}\label{closure:sampled-ray-domination}
\int F(\gamma^{x_i}-x_i)\,\mathfrak C_{n,k,T}^{\,g}(dH)\le k!K^k\,\mathsf E_{\mathsf P}F(\eta).
\end{equation}
Moreover,
\begin{equation}
\mathbb E[G_nA_n^kD_n] \le(k+1)!K^{k+1}. \label{closure:eq4_4}
\end{equation}
\end{lemma}

\begin{proof}
\emph{The bound on the insertion measure.} Conditional on $\mathcal S_k$, the remaining tree is generated by Wilson's algorithm with initial tree $\mathcal S_k$. The ball and retention restrictions on the sampled vertices are $\mathcal S_k$-measurable and preserve this conditional completion law. The marked Wilson identity gives the equality in \eqref{closure:eq4_2}. Let $y\notin\mathcal S_k$ be the next vertex, and suppose its Wilson walk first hits $\mathcal S_k$ at $w$. Assign the walk to the least index $i$ with $w\in\gamma^{x_i}$. If its loop-erased branch has length $\ell$ and $d_{\mathcal U}(0,w)+\ell\le q_n$, then
\begin{equation}\label{closure:eq4_1}
d_{\mathcal U}(x_i,y)\le d_{\mathcal U}(x_i,w)+\ell\le d_{\mathcal U}(x_i,0)+d_{\mathcal U}(0,w)+\ell\le2q_n\le\lfloor2nT\rfloor.
\end{equation}
Before hitting $w$, the same walk also avoids $\gamma^{x_i}$. Stopping it on that ray therefore gives the same loop-erased branch. Its contribution is counted by the insertion measure for $\gamma^{x_i}$ and $B_{\mathcal U}(x_i,\lfloor2nT\rfloor)$. Vertices already in $\mathcal S_k$ are assigned to their least containing ray; their atoms satisfy the same distance bound and are included in the complete ray insertion measures. Summing over the $k+1$ labels proves \eqref{closure:eq4_2}.

\emph{Moment bounds by iteration.} Put $M_{n,k}(F)=\mathbb E[G_nA_n^kF(\gamma^0)]$. After $k-1$ retained samples have been exposed, $F(\gamma^0)$ is fixed. Dropping the last sample's retention condition and applying \eqref{closure:eq4_2} to $\mathcal S_{k-1}$ gives $M_{n,k}(F)\le kK M_{n,k-1}(F)$. Since $M_{n,0}(F)\le\mathbb EF(\gamma^0)$, iteration proves \eqref{closure:eq4_3}, first for bounded $F$ and then by monotone convergence. Similarly, expose $k$ retained samples and drop the requirement that the next sample fail the retention condition. Applying \eqref{closure:eq4_2} gives $\mathbb E[G_nA_n^kD_n]\le(k+1)K \mathbb E[G_nA_n^k]$, proving \eqref{closure:eq4_4}.

\emph{Domination by the ordinary ray law.} The case $i=0$ is \eqref{closure:eq4_3}. By symmetry of the ordered samples in \eqref{closure:eq10_1}, it suffices to treat $i=1$. Apply mass transport from the original root to $x_1$, translate $x_1$ to zero, and write $y$ for the old root. The left side of \eqref{closure:sampled-ray-domination} becomes
\[
\mathbb E\!\left[g_n(0)F(\gamma^0)b_n^k\sum_{y\in B_{\mathcal U}(0,q_n)}g_n(y)\sum_{z_2,\ldots,z_k\in B_{\mathcal U}(y,q_n)}\prod_{j=2}^kg_n(z_j)\right].
\]
Expose the ray from $0$, then that from $y$, and then those from $z_2,\ldots,z_k$. The first insertion has mass at most $Y_{n,T}(\gamma^0)\le K$ on $\{g_n(0)=1\}$. All subsequent samples remain in the ball centered at $y$, and every exposed endpoint, including $0$, belongs to that ball. The proof of \eqref{closure:eq4_2}, with $y$ in place of the root, therefore gives the successive bounds $2K,\ldots,kK$. Indeed, its distance inequality uses only that all exposed endpoints are within $q_n$ of the chosen center. The ball and retention conditions are determined by the exposed rays at each step, so the conditional Wilson completion law still applies. Integrating the remaining samples in reverse order, keeping $F(\gamma^0)$ fixed, gives $k!K^k\mathbb E[g_n(0)F(\gamma^0)]$. Dropping $g_n(0)$ proves \eqref{closure:sampled-ray-domination}. Complete insertion measures include repetitions throughout; truncation and monotone convergence cover unbounded $F$.
\end{proof}

\begin{lemma}[Asymptotically inactive retention condition]
\label{closure:next-good-cutoff}
\label{the-next-goodness-cutoff-is-asymptotically-inactive}
For every fixed integer $k\ge0$, the ordered $(k+1)$-point Campbell measure with the root and first $k$ points retained and its restriction requiring the last point also to be retained differ in total variation by
\begin{equation}
\mathbb E[G_nA_n^kD_n] \le\sqrt{(2k+1)!K^{2k+1}\delta_n}\longrightarrow0. \label{closure:eq5_1}
\end{equation}
The same bound holds after recording the full marked history.
\end{lemma}

\begin{proof}
The measures are nested, and their total-mass difference is $\mathbb E[G_nA_n^kD_n]$. Cauchy-Schwarz gives $\mathbb E[G_nA_n^kD_n]\le (\mathbb E[G_nD_n])^{1/2} (\mathbb E[G_nA_n^{2k}D_n])^{1/2}$. Use $\mathbb E[G_nD_n]\le\delta_n$ and \eqref{closure:eq4_4} at order $2k$ to obtain \eqref{closure:eq5_1}. Recording the history is a measurable pushforward and cannot increase total variation.
\end{proof}

\subsection{Insertion along the sampled future rays}\label{closure:ray-reference}

Write $\mu_n=\mathfrak C_{n,k,T}^{\,g}$ for the rest of this section. A function $Z_n(H)$ tends to zero in $\mu_n$-measure if $\mu_n\{|Z_n|>a\}\to0$ for every $a>0$; it is tight if $\lim_{b\to\infty}\limsup_n\mu_n\{|Z_n|>b\}=0$. Put $x_0=0$, $R^i=\gamma^{x_i}$ and $\Gamma^i=R^i-x_i$. Thus $R^i$ is the full future ray, whereas $\gamma^i$ is the branch recorded when $x_i$ was sampled. For $i\ge1$,
\[
R^i[0,m_i]=\overleftarrow{\gamma^i},\qquad H_i:=d_{\mathcal S(H)}(0,\gamma_0^i),\qquad H_i+m_i\le q_n.
\]
Assign attachment sites to the disjoint sets $\eta[0,q_n]$ and $\gamma^i[1,m_i]$, $1\le i\le k$, so a shared attachment is counted on its earlier path. Empty branches contribute no sites. Denote their union by $\mathcal V_n=\mathcal S(H)\cap B_{\mathcal S(H)}(0,q_n)$, the assigned index by $i(x)$, and the forward coordinate of $x$ on that path divided by $n$ by $a_x$.

For a testing walk from $x$, let $J_m=|\operatorname{LE}(W[0,m])|$ and let $\Lambda_{n,m}^x$ be its stopped loop-erasure with lifetime $J_m/n$, space divided by $R_n$, and linear interpolation. Its spatial coordinates are absolute, so its starting label is $x/R_n$. The reference insertion measure is
\begin{equation}\label{closure:reference-kernel}
Q_n^{\mathrm{ray}}(H,F)=b_n\sum_{x\in\mathcal V_n}\sum_{m\ge1}\mathbb E_x\!\left[\mathbf1_{\{m<\tau_{R^{i(x)}}^+,\ d_{\mathcal S(H)}(0,x)+J_m\le q_n\}}F(i(x),a_x,\Lambda_{n,m}^x)\right].
\end{equation}
Both this measure and the rescaled complete insertion measure $\widehat Q_n(H,\cdot):=(\operatorname{Scale}_n)_*\overline Q_{n,\mathcal S(H),T}$ are measures on $\{0,\ldots,k\}\times[0,T]\times\mathscr P_T$. Equip this product with the maximum of the corresponding bounded coordinate metrics used for histories. Reference paths may hit other exposed branches, so they are treated as insertion records in this product space before comparison with Wilson insertions. For a consistent rescaled history $\mathsf h$, write $Q_T(\mathsf h,\cdot)$ for the continuum kernel \eqref{min:continuum-row} in these coordinates, taking $\ell_0=T$ and $\phi_0=\xi$.

\begin{lemma}[Convergence of the reference kernels]\label{closure:reference-convergence}
Suppose $(\operatorname{Scale}_n)_*\mu_n\Longrightarrow\mathfrak C_{k,T}$. Then
\begin{equation}\label{closure:reference-limit}
\int d_{\rm BL}\!\left(Q_n^{\mathrm{ray}}(H,\cdot),Q_T(\widehat H_n,\cdot)\right)\,\mu_n(dH)\longrightarrow0.
\end{equation}
Moreover, $\|Q_n^{\mathrm{ray}}(H,\cdot)\|\le(k+1)K$ for $\mu_n$-almost every $H$.
\end{lemma}

\begin{proof}
\emph{Restriction of the ordinary-ray kernel.} If $i\ge1$ and $x=R^i(j)$ with $0\le j<m_i$, its root distance is $H_i+m_i-j$. The length restriction in \eqref{closure:reference-kernel} implies
\[
j+J_m\le q_n-H_i-m_i+2j\le q_n-H_i+m_i\le2q_n\le\lfloor2nT\rfloor.
\]
Hence this branch's reference contribution is exactly a restricted image of $I_{n,2T}^{\Gamma^i}$ from Proposition~\ref{bridge:initial-kernel}. With $(s,(v,\beta))$ its relative insertion coordinates and $\ell_i=m_i/n$, $h_i=H_i/n$, the restriction and change of coordinates are
\begin{equation}\label{closure:reference-transform}
0\le s<\ell_i,\qquad h_i+\ell_i-s+v\le q_n/n,\qquad (s,(v,\beta))\longmapsto(i,\ell_i-s,(v,\beta_i(\ell_i-s)+\beta)).
\end{equation}
For $i=0$, restrict instead to $s+v\le q_n/n$ and translate by $\xi(s)$. Each contribution is bounded in mass by $\|I_{n,2T}^{\Gamma^i}\|\le Y_{n,2T}(\Gamma^i)\le K$, proving the asserted bound. By \eqref{closure:sampled-ray-domination} and Proposition~\ref{bridge:initial-kernel},
\[
\int d_{\rm BL}(I_{n,2T}^{\Gamma^i},I_{2T})\,d\mu_n\le k!K^k\,\mathsf E_{\mathsf P}d_{\rm BL}(I_{n,2T}^{\eta},I_{2T})\longrightarrow0.
\]

\emph{Convergence after restriction and change of coordinates.} If consistent histories $\mathsf h_n\to\mathsf h$ and finite relative insertion measures $\nu_n\Longrightarrow I_{2T}$, the restricted images in \eqref{closure:reference-transform} converge to the corresponding contribution to $Q_T(\mathsf h,\cdot)$. Indeed, under $I_{2T}$ the lifetime is $v=t$, and the possible discontinuities lie on $s=\ell_i$ or $h_i+\ell_i-s+t=T$, both of zero $ds\,dt$ measure. Uniform convergence of the attachment paths gives continuity of the changing translations and reflection away from these boundaries. Extend the change of coordinates outside the restriction by replacing $\ell_i-s$ everywhere with $(\ell_i-s)_+$. The same argument covers the initial ray and $\ell_i=0$, when the limiting contribution is zero. Changing variables $a=\ell_i-s$ gives the finite-branch limit
\[
\kappa\int_0^{\ell_i}\!da\int_0^{T-h_i-a}\!dt\;\delta_{(i,a)}\otimes\operatorname{Law}(\beta_i(a)+B^{[t]}).
\]
The same argument with $\nu_n=I_{2T}$ and root-distance bound $T$ proves continuity of $Q_T(\mathsf h,\cdot)$. Thus the distance between the restricted image and $Q_T(\mathsf h_n,\cdot)$ tends to zero along the preceding sequences. This convergence is uniform when the histories range over a compact set and $\nu_n$ approaches $I_{2T}$: otherwise a convergent subsequence of histories contradicts the preceding statement. Tightness of the scaled history measures and the displayed ordinary-ray error thus give convergence of the distance in \eqref{closure:reference-limit} in $\mu_n$-measure. The masses are bounded by $(k+1)K$ and $(k+1)\kappa T^2/2$, and $\sup_n\|\mu_n\|<\infty$, so splitting the integral at any positive threshold proves \eqref{closure:reference-limit}.
\end{proof}

\subsection{Hitting the other exposed branches}\label{closure:screening}\label{screening-by-the-rest-of-a-finite-history}

The full future ray already contains the common infinite tail of the exposed tree. The rest of the exposed tree therefore consists only of finite branches. We compare the two insertion measures by omitting attachment sites near recorded endpoints and using spatial separation to bound hits on the other branches.

\begin{lemma}[Hitting a spatially separated finite set]\label{closure:separated-hitting}
Let $W$ start at $x$, let $\mathcal C_x$ be a finite set containing $x$, and put $E_x=\{W[1,r_n]\cap\mathcal C_x=\varnothing\}$ and $p_x=\mathbb P_x(E_x)$. If $A_*$ is finite and $\operatorname{dist}(x,A_*)\ge\rho R_n$, then, for every $m\ge r_n$,
\begin{equation}\label{closure:eq9_1}
\mathbb P_x(E_x,\ W[1,m]\cap A_*\ne\varnothing)\le Cp_x\frac{\operatorname{cap}(A_*)}{\rho^2R_n^2}+C\exp(-\rho^2R_n^2/(Cr_n)).
\end{equation}
\end{lemma}
\begin{proof}
By \eqref{bridge:gaussian-maximal}, the first $r_n$ steps stay within $\rho R_n/2$ of $x$ except with the probability in the second term. If those steps stay within that ball, $W$ has not hit $A_*$ and $\operatorname{dist}(W_{r_n},A_*)\ge\rho R_n/2$. The Markov property and the Green-function hitting bound \cite[equations~(2.1)--(2.2)]{HutchcroftSousi2023} bound the subsequent hitting probability by $C\operatorname{cap}(A_*)/(\rho^2R_n^2)$. Integrating over $E_x$ proves the result.
\end{proof}

\begin{lemma}[Separation under convergence of history coordinates]\label{closure:history-geometric-separation}
Fix $A>T$. For a consistent history $\mathsf h$ truncated at $A$, let $\mathcal T_A(\mathsf h)$ be its finite tree, $I_i(\mathsf h)$ its recorded arcs, and $J(\mathsf h)$ the set of all arc endpoints, including attachments and zero-length arcs. Write $\phi_{\mathsf h}$ for its spatial label map and, for $\delta>0$, set
\[
C_{i,\delta}(\mathsf h)=\{x\in I_i(\mathsf h):d_{\mathcal T_A(\mathsf h)}(x,J(\mathsf h))\ge\delta\},\qquad \operatorname{Sep}_{\delta}(\mathsf h)=\inf_{\substack{i\ne j,\ x\in C_{i,\delta}(\mathsf h)\\y\in I_j(\mathsf h)}}|\phi_{\mathsf h}(x)-\phi_{\mathsf h}(y)|,
\]
with $\inf\varnothing=+\infty$. If $\mathsf h_n\to\mathsf h$ and $\phi_{\mathsf h}$ is injective, then $\liminf_n\operatorname{Sep}_{\delta}(\mathsf h_n)\ge\operatorname{Sep}_{\delta/2}(\mathsf h)>0$.
\end{lemma}
\begin{proof}
The parent indices are eventually fixed. Parameterize every arc on $[0,1]$, using a constant parameterization for zero-length arcs. Recursive gluing expresses all tree distances continuously in the finitely many lengths and attachment coordinates, uniformly over these parameters. Distances to $J$ therefore converge uniformly too: every moving attachment is an endpoint of its child arc. The spatial labels also converge uniformly. Hence a point at distance at least $\delta$ from $J(\mathsf h_n)$ corresponds, for large $n$, to a limiting point at distance at least $\delta/2$ from $J(\mathsf h)$, and its distance from any competing spatial arc differs by $o(1)$. These retained limiting points are distinct from all points of the competing arcs. Injectivity and compactness give positive separation, proving the inequality. If there are no such limiting pairs, there are eventually none before the limit either.
\end{proof}

\begin{proposition}[Comparison with the full exposed tree]\label{closure:history-screening}
Suppose $(\operatorname{Scale}_n)_*\mu_n\Longrightarrow\mathfrak C_{k,T}$. Then
\begin{equation}\label{closure:full-ray-comparison}
\int\|\widehat Q_n(H,\cdot)-Q_n^{\mathrm{ray}}(H,\cdot)\|_{\rm TV}\,\mu_n(dH)\longrightarrow0.
\end{equation}
\end{proposition}
\begin{proof}
\emph{Reduction to hitting a finite set.} Put $\mathcal A_i=\mathcal S(H)\setminus R^i$. Every added branch lies in the original root ball, and $R^i$ meets the initial ray at root height at most $q_n$. It therefore contains that ray beyond height $q_n$, so
\begin{equation}\label{closure:finite-competitor}
\mathcal A_i\subseteq\mathcal V_n,\qquad |\mathcal A_i|\le(k+1)(q_n+1),\qquad \operatorname{cap}(\mathcal A_i)\le C_{k,T}n.
\end{equation}
Before recording the loop-erasure, the positive-length full-tree measure is dominated by the reference measure: their difference counts walks that avoid $R^{i(x)}$ but hit $\mathcal A_{i(x)}$. Denote this nonnegative occupation mass, with the actual length restriction retained, by $\Delta_n(H)$. Recording the same attachment and path coordinates cannot increase total variation, and the exposed-vertex atoms have mass at most $b_n|\mathcal V_n|=o(1)$. It suffices to show $\Delta_n\to0$ in $\mu_n$-measure.

\emph{Local avoidance weights and time cutoffs.} For $0\le j\le q_n$ on each full ray, put $\mathcal C_{i,j}=R^i[(j-c_n)_+,j+c_n]$ and $p_{i,j}=\mathbb P_{R^i(j)}\{W[1,r_n]\cap\mathcal C_{i,j}=\varnothing\}$. Proposition~\ref{bridge:row-convergence} has diffuse limit $\kappa\,ds$, so its cumulative masses converge uniformly on $[0,T]$. Transferring the corresponding error event by \eqref{closure:sampled-ray-domination} gives
\begin{equation}\label{closure:uniform-cumulative-history-rows}
\omega_n:=\max_{0\le i\le k}\sup_{v\in[0,T]}\left|\frac{L_n^2}{n}\sum_{j=0}^{\lfloor nv\rfloor}p_{i,j}-\kappa v\right|\longrightarrow0\quad\text{in }\mu_n\text{-measure}.
\end{equation}
For an assigned site $x=R^{i(x)}(j)$, use this ray's window and denote its early avoidance probability by $p_x$. Then $E_n=(L_n^2/n)\sum_{x\in\mathcal V_n}p_x$ is tight, since it is bounded by the sum of the $k+1$ full-ray cumulative masses. There are $O_{k,T}(n)$ sites and every inserted length is at most $q_n$. Lemma~\ref{bridge:raw-cutoffs}, with killing sets $R^{i(x)}$ and their windows, gives for a sufficiently large fixed $M$ a total reference contribution outside $r_n\le m\le MN_n$ bounded by
\[
C_{k,T}\left(L_n^{-11}+\frac{\log\log n}{L_n}E_n+L_n^{-1}\right)=o_{\mu_n}(1).
\]
Here and below $o_{\mu_n}(1)$ denotes convergence to zero in $\mu_n$-measure. The short-time term also follows directly by counting at most $r_n$ times per site.

\emph{Sites near recorded endpoints.} Fix $A=T+1$ and $\delta>0$, and truncate the initial ray at $A$ in rescaled coordinates. Omit assigned sites within intrinsic distance $\delta$ of any recorded endpoint in this finite tree. On each recorded arc these neighborhoods are a bounded number of intervals of total rescaled length at most $C_k\delta$. Reflection into the corresponding full-ray coordinate preserves these lengths. Equation~\eqref{closure:uniform-cumulative-history-rows} bounds their normalized early-avoidance mass by $C_k\delta+o_{\mu_n}(1)$. For $m\ge r_n$, reference avoidance implies early avoidance, so their contribution up to $MN_n$ is at most $C_{k,T,M}\delta+o_{\mu_n}(1)$.

\emph{Spatial separation and hitting estimates.} Under the continuum Campbell measure, conditional on the genealogy and lengths, the edge increments are independent Brownian motions. Each geodesic between recorded endpoints, read relative to its first label, is therefore four-dimensional Brownian motion, including any reversed traversals, and is injective almost surely \cite[Theorem~2, p.~80]{DvoretzkyErdosKakutani1950}. Finitely many such geodesics cover every pair of tree points. Thus the finite limiting tree has injective labels, after zero-length arcs are contracted. Lemma~\ref{closure:history-geometric-separation} and weak convergence now imply that, for every $\theta>0$, one can choose $\rho>0$ such that
\[
\limsup_n\mu_n\{\operatorname{Sep}_{\delta}(\mathsf R_A\widehat H_n)\le\rho\}<\theta.
\]
For completeness, choose a compact set of injectively labelled limiting histories with complement of mass less than $\theta$ and a positive lower bound for $\operatorname{Sep}_{\delta/2}$ there. The preceding lemma supplies a neighborhood with a uniform positive lower bound for $\operatorname{Sep}_{\delta}$; Portmanteau gives the displayed assertion.

Outside this exceptional set, every remaining site $x$ is at least $\rho R_n$ from $\mathcal A_{i(x)}$, since $R^{i(x)}$ contains its entire assigned arc. Apply Lemma~\ref{closure:separated-hitting} with its full-ray window and \eqref{closure:finite-competitor}. Summing over sites and $r_n\le m\le MN_n$, using $b_nN_n=L_n^2/n$, $R_n^2=nL_n/\alpha$ and $R_n^2/r_n\asymp L_n^{13}$, gives
\[
\Delta_n\le C_{k,T,M}\delta+\frac{C_{k,T,M}}{\rho^2L_n}E_n+C_{k,T,M}L_n^2e^{-c\rho^2L_n^{13}}+o_{\mu_n}(1).
\]
Given a positive error threshold, first choose $\delta$ small, then $\rho$ so that the exceptional history mass is below $\theta$, and finally let $n\to\infty$. Tightness of $E_n$ makes the second term vanish. Letting $\theta\downarrow0$ proves $\Delta_n\to0$ in $\mu_n$-measure. Both insertion masses are at most $(k+1)K$ by Lemmas~\ref{closure:bounded-rows} and \ref{closure:reference-convergence}. Since $\sup_n\|\mu_n\|<\infty$, convergence in measure implies \eqref{closure:full-ray-comparison}.
\end{proof}

\subsection{One-step insertion and proof of Campbell convergence}\label{closure:induction}\label{the-finite-history-induction}

We now prove Proposition~\ref{closure:campbell-convergence} by combining the exact Campbell recursion with a one-step insertion limit. For a history measure $\mu$ and insertion kernel $Q$, write $\mu\ltimes Q$ for the pushforward of their joint measure under the map that appends the insertion record to its history.

\begin{lemma}[Exact discrete Campbell recursion]
\label{closure:discrete-recursion}
Let $\mathfrak C_{n,k+1,T}^{\,g,-}$ be the history measure in \eqref{closure:eq10_1} with the root and first $k$ samples retained, and no retention condition on the last sample. Keeping the attachment vertex and lattice path in the Wilson insertion measure, one has
\begin{equation}\label{closure:exact-campbell-recursion}
\mathfrak C_{n,k+1,T}^{\,g,-} =\mathfrak C_{n,k,T}^{\,g}\ltimes\overline Q_{n,\mathcal S(H),T}.
\end{equation}
Retaining the last endpoint gives precisely the restricted measure $\mathfrak C_{n,k+1,T}^{\,g}$, and
\[
\|\mathfrak C_{n,k+1,T}^{\,g,-}-\mathfrak C_{n,k+1,T}^{\,g}\|_{\rm TV} =\mathbb E[G_nA_n^kD_n]\longrightarrow0.
\]
Applying $(\operatorname{Scale}_n)_*$ gives the corresponding identity for rescaled histories.
\end{lemma}

\begin{proof}
For a fixed tuple $(0,x_1,\ldots,x_k)$, expose its future rays by Wilson's algorithm. Their union $\mathcal S$ determines the root-to-sample paths and all recorded future rays. Hence the ball and retention restrictions are $\mathcal S$-measurable: conditioning or reweighting by them preserves the conditional Wilson completion law. The marked occupation identity supplies its complete next insertion measure, including exposed-vertex atoms. Summing over tuples and applying Tonelli gives \eqref{closure:exact-campbell-recursion}. The appended history determines the new endpoint's future ray and whether it is retained. Lemma~\ref{closure:next-good-cutoff} therefore gives the displayed total-variation bound.
\end{proof}

\begin{proposition}[One insertion under a retained history measure]\label{closure:one-step}
Fix $T,K$ as above and $k\ge0$, and suppose $(\operatorname{Scale}_n)_*\mu_n\Longrightarrow\mathfrak C_{k,T}$ for $\mu_n=\mathfrak C_{n,k,T}^{\,g}$. Then
\begin{equation}\label{closure:one-step-kernel}
\int d_{\rm BL}\!\left(\widehat Q_n(H,\cdot),Q_T(\widehat H_n,\cdot)\right)\,\mu_n(dH)\longrightarrow0,
\end{equation}
and, writing $\overline Q_n(H,\cdot)=\overline Q_{n,\mathcal S(H),T}$,
\begin{equation}\label{closure:one-step-history}
(\operatorname{Scale}_n)_*(\mu_n\ltimes\overline Q_n)\Longrightarrow\mathfrak C_{k,T}\ltimes Q_T=\mathfrak C_{k+1,T}.
\end{equation}
\end{proposition}
\begin{proof}
The triangle inequality, \eqref{closure:reference-limit} and \eqref{closure:full-ray-comparison} give \eqref{closure:one-step-kernel}. To append the new coordinates, let $F$ be a unit bounded-Lipschitz test on $\mathscr X_{k+1,T}$, and denote by $F_{\mathsf h}(i,s,(t,\beta))$ its value after this insertion is appended to $\mathsf h$. The maximum product metric makes $F_{\mathsf h}$ a unit bounded-Lipschitz test in the insertion coordinates, uniformly in $\mathsf h$. Thus \eqref{closure:one-step-kernel} yields
\[
(\operatorname{Scale}_n)_*(\mu_n\ltimes\overline Q_n)(F)=\int G_F(\widehat H_n)\,\mu_n(dH)+o(1),
\]
where $G_F(\mathsf h)=Q_T(\mathsf h,F_{\mathsf h})$ on consistent histories. Extend it to the whole coordinate space by
\[
G_F(\mathsf h)=\kappa\sum_{i=0}^k\int_{[0,T]^2}\!\mathbf1_{\{s\le\ell_i,\ h_i+s+t\le T\}}\mathbb E F_{\mathsf h}\!\left(i,s,(t,\phi_i(s)+B_{\cdot\wedge t})\right)\,ds\,dt.
\]
For convergent history coordinates, the parent indices are eventually fixed, the heights and lengths converge, and the attachment paths converge uniformly on $[0,T]$. Coupling with the same Brownian motion makes the integrand converge except on the two zero-area boundaries $s=\ell_i$ and $h_i+s+t=T$. Dominated convergence proves that $G_F$ is continuous, including at zero-length branches, and $|G_F|\le\kappa(k+1)T^2/2$. Weak convergence of the scaled histories and the continuum recursion in Corollary~\ref{min:continuum-campbell} therefore give
\[
\int G_F(\widehat H_n)\,\mu_n(dH)\longrightarrow\mathfrak C_{k,T}(G_F)=\mathfrak C_{k+1,T}(F),
\]
which proves \eqref{closure:one-step-history}.
\end{proof}

\begin{proof}[Proof of Proposition~\ref{closure:campbell-convergence}]
We prove \eqref{closure:campbell-limit} by induction on $k$. At $k=0$, the initial-ray measure is $G_n\mathsf P$, and its total-variation distance from $\mathsf P$ is $\mathsf P(G_n=0)\to0$ by Lemma~\ref{closure:good-deletion}. The initial-ray invariance principle from Section~\ref{prelim:lerw} gives the required limit $\mathfrak C_{0,T}$.

Suppose the convergence holds at order $k$. Proposition~\ref{closure:one-step} gives convergence after one complete insertion. Lemma~\ref{closure:discrete-recursion} identifies this history measure with $\mathfrak C_{n,k+1,T}^{\,g,-}$, and retaining the new endpoint changes it by $o(1)$ in total variation, by Lemma~\ref{closure:next-good-cutoff}. Its scaled image therefore converges to $\mathfrak C_{k+1,T}$, completing the induction. The sampled-ray domination used at each step follows from Lemma~\ref{closure:bounded-rows} independently of this convergence.

Finally, the root-to-sample distances, sampled distance arrays and spatial endpoint labels are continuous functions of the history coordinates: recursive gluing gives continuous tree distances, as in Lemma~\ref{closure:history-geometric-separation}, and endpoint evaluation is continuous for stopped paths. Zero-length branches cause no exception. For $k\ge1$, $G_n^k=G_n$, so projecting $\mathfrak C_{n,k,T}^{\,g}$ onto these arrays gives exactly the expected $k$th rooted marked moment measure of the retained ball. The last assertion follows by continuous pushforward.
\end{proof}

The induction has identified all expected finite-sample measures of the retained balls. We now recover each measured ball's law from its moments, restore the removed mass, and then let the radius increase.

\begin{proof}[Proof of Theorem~\ref{min:main}]
Fix $T>0$, and let $\widetilde{\mathscr U}_{n,T}$ be the rescaled intrinsic ball with retained measure $\widetilde\mu_{n,T}$. Proposition~\ref{closure:campbell-convergence} gives convergence of all its positive-order rooted marked moment measures. The limiting ball volume satisfies \eqref{min:continuum-factorial}, so Lemma~\ref{prelim:marked-moments} yields
\[
\widetilde{\mathscr U}_{n,T} \Longrightarrow B_T(\mathcal K_{\kappa,I_4/4}).
\]
Write $\mathscr U_{n,T}$ for the same ball with its original vertex measure. Embed both balls into the same metric space, keeping their spatial labels. Their roots agree, and their measures differ by a nonnegative measure of total mass $V_n-G_nA_n$. The common embedding and Lemma~\ref{closure:good-deletion} therefore give
\begin{equation}
d_{\mathrm{r,mGP}}(\mathscr U_{n,T},\widetilde{\mathscr U}_{n,T}) \le V_n-G_nA_n\longrightarrow0\quad\text{in }L^1. \label{closure:measure-deletion-distance}
\end{equation}
Slutsky's theorem therefore gives $\mathscr U_{n,T}\Longrightarrow B_T(\mathcal K_{\kappa,I_4/4})$.

To pass from balls to the full space, use the continuous identity from finite rooted marked spaces to the local state space, see the paragraph ``Two facts about the local topology'' in Section~\ref{prelim:continuum}. Thus the fixed-$T$ convergence also holds in the local topology. Restriction at radius $T$ changes none of the smaller measured balls, so deterministically
\[
d_{\mathrm{loc}}(\mathscr U_n,\mathscr U_{n,T})\le e^{-T},\qquad d_{\mathrm{loc}}\bigl(\mathcal K_{\kappa,I_4/4}, B_T(\mathcal K_{\kappa,I_4/4})\bigr)\le e^{-T}.
\]
A converging-together argument, first letting $n\to\infty$ at fixed $T$ and then $T\to\infty$, proves the asserted local convergence.
\end{proof}

\begin{proof}[Proof of Corollary~\ref{min:volume-occupation-law}]
Projecting the one-point measure \eqref{min:distance-matrix} onto the distance from the root and the spatial label, the fixed-radius convergence proved above gives \eqref{min:occupation-limit}, with the limit identified by Proposition~\ref{prelim:icsbm-volume}. Continuity of total mass gives the volume convergence.
\end{proof}

\bigskip
\setlength{\bibindent}{20pt}

\end{document}